\documentclass[a4paper, 12pt]{article}

\usepackage[dvips]{graphics,color}
\def\spc{\color[rgb]{1,0.2,0.2}}
\def\nc{\normalcolor}

\usepackage{amsmath}
\usepackage{amssymb}
\usepackage{graphicx}

\usepackage{verbatim}
\usepackage[british]{babel}

\newtheorem{lemma}{Lemma}

\newtheorem{theorem}{Theorem}
\newtheorem{rem}{Remark}

\newcommand{\R}{{\mathbb{R}}}
\renewcommand{\P}{\mathsf{P}}
\newcommand{\Z}{\mathbb{Z}}
\newcommand{\eps}{\varepsilon}
\newcommand{\E}{\mathsf{E}}
\newcommand{\G}{\mathsf{L}}
\renewcommand{\L}{\Lambda}

\def\bn#1\en{\begin{align*}#1\end{align*}}
\def\bnn#1\enn{\begin{align}#1\end{align}}

\def\|{\, |\, }
\begin{document}

\title{Long term  behaviour  of 
 locally interacting  birth-and-death processes}

\author{
Vadim Shcherbakov\footnote{ Royal Holloway, University of London. Email: Vadim.Shcherbakov@rhul.ac.uk
}
\, and 
Stanislav Volkov\footnote{
Lund University.  Email: S.Volkov@maths.lth.se
}
}

\maketitle

\begin{abstract}
In this  paper we study long-term evolution of a finite 
system of locally interacting birth-and-death processes labelled by vertices of a finite connected graph. A detailed description of the asymptotic behaviour is obtained in the case of both constant vertex degree graphs and star graphs. The model is motivated by modelling interactions between  populations and is related to interacting particle systems,  Gibbs models with unbounded spins, as well as urn models with interaction.
\end{abstract}

\section{The model}\label{intro}
 Let  $\Lambda$ be  a finite connected graph.  If two vertices $x, y\in \L$ are connected by an edge, call them {\it neighbours} and write $x\sim y$.  Let  $\Z$ be the set of all  integers  and  $\Z_{+}$ be the set of all non-negative integers including zero. Consider a continuous time Markov chain (CTMC)  $\xi(t)=\{\xi_x(t),\, x\in \Lambda\}\in \Z_{+}^{\Lambda}$ with the following transition  rates:  given $\xi(t)=\xi\in \Z_{+}^{\Lambda}$ a component (a spin)  $\xi_x$ 
increases by $1$ at the rate $e^{\alpha\xi_x+\beta\phi(x,\xi)}$, where  $\alpha, \beta\in \R$,
\begin{equation}
\label{phi}
\phi(x, \xi)=\sum\limits_{y:y\sim x}\xi_y
\end{equation}
and at the same time each positive component  $\xi_x$ decreases by $1$ at constant rate $1$. 

This birth-and-death  dynamics belongs to a class   of  stochastic dynamics  which is 
 used in  statistical physics  to  describe  the time evolution of a system of interacting spins.
Our particular dynamics is motivated by adsorption-desorption processes, where adsorption rates depend
on a local environment and  an adsorbed particle can depart at a non-zero rate (\cite{Evans}).
 It is  closely related to a particle deposition  on a discrete substrate  and  urn models with interaction
(e.g.,  \cite{Kulske}, \cite{Volk1}, and~\cite{Volk2}).
Recall also  that a  birth-and-death process on the non-negative integer half-line  is a classic probabilistic  model for the population size so that 
the  Markov chain can be  used for modelling different types of  interaction between populations,
where  a component $\xi_x(t)$ can be interpreted  as the size of a population which is located at $x\in \Lambda$ at time~$t$.

 If we assume that spins are bounded and consider the same birth-and-death dynamics then 
we will get   a finite  ergodic Markov chain  whose equilibrium distribution  is a Gibbs measure (see Remark~\ref{Gibbs}). 
A particular case  of the model with bounded spins, where   $\alpha=\beta,\, \Lambda\subset \Z^d$,   was studied in~\cite{Yamb}.
For instance, if a spin  takes   values $0$ and $1$ only, and, in addition, $\alpha=\beta>0$, then 
we obtain a finite Markov chain whose equilibrium distribution is a Gibbs measure on  $\{0, 1\}^\L$ which is equivalent 
to a particular case of the famous Ising model on $\{-1, 1\}^{\L}$.    
The  main goal in~\cite{Yamb} was to  study the asymptotic  behaviour of the stationary distribution 
  as~$\L\uparrow \Z^d$.  In general, the   asymptotic behaviour  of such equilibrium distributions   in 
thermodynamic limit, i.e.  as  graph $\L$ expands,   is  of  main interest in statistical physics. 

The aim of this paper, on the other hand, is to describe the asymptotic behaviour of  the  Markov chain with {\it unbounded} spins
 as time tends to infinity while the underlying graph remains fixed.
In this case  we deal with a countable Markov chain 
that can be either  recurrent (or even ergodic)  or non-recurrent (e.g., transient, or even explosive) depending both on 
the  graph $\L$ and the values of  parameters~$\alpha, \beta$.

It is easy to see that if $\beta=0$ then the structure of graph $\L$ is irrelevant and 
the components of CTMC $\xi(t)$ are independent identically distributed birth-and-death processes with values  in $\Z_{+}$. 
The well known  results  for birth-and-death processes (e.g.\ see~\cite{Karlin} or~\cite{Ross}) yield that
  if $\alpha>0, \beta=0$, then  each component  is an explosive MC. In turn, it  implies  that CTMC $\xi(t)$ is  explosive. 
Moreover,  independence of spins imply that their times to explosion are also independent
and this allows to  repeat   the well  known  Rubin's argument 
(used in~\cite{Davis}  in the case of classic Polya urn scheme) in order to obtain 
 that with probability $1$ only a single component of $\xi(t)$ explodes. 
Notice that this fact can be also inferred from our Theorem~\ref{T2}.  
A non-zero interaction does not change   the explosive behaviour of 
the Markov chain in the case $\alpha>0$ but  escape 
to infinity can happen   in various  ways  which    depend on  both  $\beta$ and  $\L$.
 
If  $\alpha<0, \beta=0$, then  CMTC $\xi(t)$ is formed by a collection of  independent  ergodic Markov chains.
It is quite obvious that if both $\alpha<0$ and  $\beta<0$ then  the Markov chain remains to be  ergodic.
If $\beta>0$, then one could intuitively expect that given $\alpha<0$ there exists some critical value $\beta_{cr}$ such that if $\beta<\beta_{cr}$, then 
 the stable ergodic evolution of the system  is still observed,  and,  in contrast,  if $\beta>\beta_{cr}$, then 
  the system becomes unstable, i.e. transient or even  explosive.
We compute this critical value explicitly in some cases.
It turns out that  $\beta_{cr}=-\alpha c(\L)$, where  $c(\L)=\nu^{-1}$ in the case of a  graph $\L$ with the constant vertex degree $\nu$
and $c(\L)=n^{-\frac{1}{2}}$ in the case of a star graph $\L$ with $n+1$ vertices.


The Markov chain under consideration  is reversible, therefore the computation of its invariant measure is straightforward.
Stationary probability  distributions arising in positive recurrent cases are  Gibbs measures with unbounded positive spins    on a  finite graph with empty boundary conditions. Consequently the model in positive recurrent cases is  closely related to  Gibbs random fields 
  with unbounded spins on graphs (see~\cite{Pasurek}, \cite{Lebowitz}, and references therein).

We  give a   detailed description of how the Markov chain escapes to infinity in all the transient cases that we consider.  We show  that due to a rapid increase of birth rates in explosive  cases, there are no death events in the system after some finite random moment of time, and the dynamics of the Markov chain is that of a pure birth process, obtained by setting the death rates to zero. 

We  will start with  results that are valid in the case of an  arbitrary finite connected  graph~$\L$; they are presented in Theorems~\ref{T1},~\ref{T2} and~\ref{T.no.triangle}. 
We also study  two special cases  in more detail, 
namely constant vertex degree graphs and  star graphs. The results for these two  cases are found in Theorems~\ref{T3},~\ref{T31} and~\ref{T4}.
Graphs with the constant vertex degree and star graphs are particular examples of spatially homogeneous graphs and of spatially inhomogeneous graphs, respectively. 
Despite the obvious difference in the structure of these graphs the long term behaviour of the corresponding Markov chains 
is similar to each other. 
 The main  features of the model dynamics are illustrated  in Section~\ref{exp}  by  a model with  graph~$\L$  formed by just two neighbouring vertices.
Proofs are given in Section~\ref{proofs}.

Finally, we denote by   $C_i,\, i=1,2,...,$ or just $C$   various constants whose exact values  are immaterial and can change from line to line.

\section{Results}
\label{results}
Let $\L$ be an arbitrary graph. Given  $\xi\in \Z_{+}^{\L}$ define  {\it potential} $U(x, \xi)$ of a vertex $x\in \L$ as the following quantity  
\begin{equation}\label{U}
U(x, \xi)=\alpha\xi_x+\beta\phi(x, \xi).
\end{equation}
Notice the  following identity  
\begin{equation}\label{Usum}
\sum_{x\in \L}U(x, \xi)=\sum_{x\in\L}(\alpha+\beta\nu(x))\xi_x,
\end{equation}
where $\nu(x)$ is  the degree  of vertex $x\in \L$, i.e.\ the number of  edges incident to the vertex.
Throughout the paper we will also denote  by  $1_{A}$  the indicator of  a set (or  event)  $A$.
In these notations, given $\xi(t)=\xi\in\Z_{+}^{\L}$ a  component $\xi_x$ jumps up by $1$ with intensity $e^{U(x, \xi)}$ and the generator of the Markov chain is therefore
\begin{equation}
\label{generator}
\G f(\xi)=\sum_{x\in \L} \left(f\left(\xi+{\sf e}^{\left(x\right)}\right)-f(\xi)\right)e^{U(x, \xi)}
+\left(f\left(\xi-{\sf e}^{\left(x\right)}\right)-f(\xi)\right)1_{\{\xi_x>0\}},
\end{equation}
where ${\sf e}^{\left(x\right)}$ is a configuration such that ${\sf e}^{\left(x\right)}_x=1$ and ${\sf e}^{\left(x\right)}_y=0$ for all $y\neq x$ (addition of configurations is understood component-wise).

Recall that the embedded Markov chain, corresponding to a continuous time Markov chain (CTMC), is a discrete time Markov chain (DTMC) with the same state space, and that makes the same jumps as the continuous time Markov chain  with probabilities proportional to the corresponding jump rates. Let~$\zeta(t)$ be the DTMC corresponding to CTMC~$\xi(t)$. The states of the embedded Markov chain will be denoted by~$\zeta$ and we will use the same symbol~$t=0, 1,2,\dots,$ to denote the discrete time.

It is easy to see that if $\beta=0$ then the components of CTMC $\xi(t)$ are independent identically distributed birth-and-death processes  in $\Z_{+}$.  It is easy to see that CTMC $\xi(t)$ is  ergodic, if $\alpha<0$, and is explosive if $\alpha>0$ respectively. Also, if both $\alpha=0$ and $\beta=0$ then CTMC $\xi(t)$   is formed by a collection of independent  simple  symmetric random walks in $\Z_{+}$ reflected at the origin. This CTMC is null recurrent if $|\L|=1,2$, and is transient if $|\L|\geq 3$. We exclude these trivial cases in what follows.  

Let us define the following function
\begin{equation}
\label{W}
W(\xi)=\exp\left(\alpha\sum_{x}\xi_x(\xi_x-1)/2+\beta\sum_{x\sim y}\xi_x\xi_y\right),\quad
\xi\in \Z_{+}^{\L}.
\end{equation}
It is easy to see that 
$$e^{U(x, \xi)}e^{W(\xi)}=e^{W\left(\xi+{\sf e}^{\left(x\right)}\right)}$$
for all $x\in \L$ and $\xi\in \Z_{+}^{\L}$.
This equation is  a detailed balance condition which implies that the Markov chain is time-reversible 
with  invariant measure $e^{W(\xi)},\, \xi\in \Z_{+}^{\L}$.
According to e.g. Theorem~1.2.4  in~\cite{FMM}, an irreducible countable Markov chain is positive recurrent (i.e. ergodic) if and only if there exists a  stationary probability  distribution, and if the latter exists then the distribution of the Markov chain converges to it as time goes to infinity. Therefore, if 
\begin{equation}
\label{Z}
Z_{\alpha, \beta, \L}=\sum\limits_{\xi\in  \Z_{+}^{\L}}e^{W(\xi)}<\infty,
\end{equation}
then  CTMC $\xi(t)$  is  ergodic   
with the stationary probability distribution given by
\begin{equation}\label{inv}
\mu_{\alpha, \beta, \L}(\xi)=\frac{e^{W(\xi)}}{Z_{\alpha, \beta, \L}} ,\quad
\xi\in \Z_{+}^{\L}.
\end{equation}
\begin{rem}
\label{Gibbs}
{\rm
If a component of the Markov chain takes values in $\{0, 1, \ldots, N\}$, where $N\geq 1$, then 
 the invariant probability distribution of the Markov chain is defined similar to measure~(\ref{inv}). Namely,
it  is a probability measure on $\{0, 1, \ldots, N\}^{\L}$ that is equal, up to a normalizing constant, 
 to function $e^{W(\xi)}$, where, in turn, 
 function $W$ is defined, as before,  by~(\ref{W}).
}
\end{rem}
Denote
\begin{align}
\label{Q}
Q(\xi)&=-\alpha\sum_{x}\xi^2_x-2\beta\sum_{x\sim y}\xi_x\xi_y,\\
\label{S}
S(\xi)&=\sum\limits_{x}\xi_x.
\end{align}
Then we can rewrite function~(\ref{W}) as
\begin{equation}
\label{QS}
W(\xi)=-\frac{1}{2}(Q(\xi) + \alpha S(\xi)).
\end{equation}
Recall that  $\nu(x)$ denotes  the degree  of vertex $x\in \L$ and notice  the following useful representations of the quadratic part 
of $W$  
\begin{align}
\label{Qform1}
Q(\xi)&=\sum_{x}(-\alpha-\beta\nu(x))\xi^2_x +
\beta\sum\limits_{x\sim y}(\xi_x-\xi_y)^2\\
&=\sum\limits_{x\in \L}\bigl(-\alpha\xi_x^2-\beta\xi_x\phi(x,\xi)\bigr)=-\sum\limits_{x\in \L}\xi_xU(x, \xi).\nonumber
\end{align}

We are ready now to formulate the findings of our paper. We start  with the results that are valid for all finite connected graphs. 
\begin{theorem}
\label{T1}
 Let  $\L$ be  a finite  connected graph.
\begin{itemize}
\item[1)] 
 If  $\alpha<0$  and $\alpha+\beta\max_{x\in \L}\nu(x)\leq 0$  then  CTMC $\xi(t)$ is not explosive. Moreover, if $\alpha<0$  and $\alpha+\beta\max_{x\in \L}\nu(x)< 0$  then CTMC $\xi(t)$ is ergodic.
\item[2)] If $\alpha\geq 0$ then CTMC $\xi(t)$ is not  ergodic.
\end{itemize}
\end{theorem}
The transient behaviour of the Markov chain in Part 2) of Theorem~\ref{T1} can 
be described more precisely under certain additional assumptions. In order to do so, define the  following event related to DTMC~$\zeta(t)$:
\begin{align}
\label{B}
B=\left\{\exists\mbox{$\tau\in\Z_+$ and a vertex $x\in \L$ such that } \right.\nonumber
\\
\left. 
\zeta_{y}(\tau+s+1)=\zeta_y(\tau+s)+1_{\{y=x\}},\,\, \forall s\geq 1\right\},
\end{align}
in other words, the process grows only at point $x$ after time $\tau$.

\begin{theorem}
\label{T2}
Let  $\L$ be  a finite graph. If $\alpha>\max\{0,\beta\}$ then 
with probability $1$ event $B$ defined by~(\ref{B}) occurs, and a single component of  CTMC $\xi(t)$ explodes.
\end{theorem}
\begin{rem}
{\rm 
Notice that we do not assume connectedness of the underlying graph in Theorem~\ref{T2}. 
}
\end{rem}
Furthermore, given $x_1, x_2\in \L$ define the following event
\begin{equation*}
B_{x_1, x_2}=\left\{\exists s\in\Z_+: \zeta_y(t)=\zeta_y(s)\text{ for all }y\notin\{x_1,x_2\}
\text{ and all } t\geq s;\right.
\end{equation*}
\begin{equation}\label{eq.just.two1}
\left.
 \lim_{t\to\infty} \frac{\zeta_{x_1}(t)}{t}=\lim_{t\to\infty}
 \frac{\zeta_{x_2}(t)}{t}=\frac{1}{2}\right\}.
\end{equation}
\begin{theorem}
\label{T.no.triangle}
Let  $\L$ be  a
 finite  connected graph without triangles, i.e.\ such that there are no three distinct vertices $x,y,z\in\L$ such that $x\sim y$, $y\sim z$ and $z\sim x$. If $0<\alpha<\beta$ then with probability $1$ there are two adjacent vertices $x_1$ and $x_2$ such that the event~(\ref{eq.just.two1}) occurs. This implies that with probability $1$ only a pair of adjacent components of the CTMC explodes.
\end{theorem}

\begin{theorem}
\label{T3}
 Let  $\L$ be  a graph with the constant vertex degree $\nu(x)\equiv \nu$.  
 \begin{itemize}
 \item[1)] 
 CTMC $\xi(t)$ is ergodic if and only if $\alpha<0$ and $\alpha+\beta\nu<0$.
 \item[2)] If $\alpha<0$ and $\alpha+\beta\nu=0$ then CTMC $\xi(t)$ is transient.
 \item[3)] If   $\alpha<0$ and $\alpha+\beta\nu>0$ then  CTMC  $\xi(t)$ is explosive.
 \item[4)] If   $\alpha>0$ then  CTMC  $\xi(t)$ is explosive. Moreover, 
\begin{itemize}
 \item[i)] if $\beta<\alpha$ then with probability $1$ the event~(\ref{B}) occurs  and a single component of CTMC $\xi(t)$ explodes;
 \item[ii)] if  $\alpha<\beta$ and the graph $\L$ is without triangles (as explained in Theorem~\ref{T.no.triangle}) then with probability $1$ the event $B_{x_1, x_2}$ occurs for some adjacent vertices $x_1, x_2\in \L$, so that with probability $1$ a pair of adjacent components  of the CTMC explodes.
\end{itemize}
\end{itemize}
\end{theorem}
Let us mention two  examples of   constant vertex degree  graphs, both with and without triangles. 
\begin{itemize}
\item[ a)] {\it  Lattice models with  local interaction.}  Let $\Z$   be  the set of all  integers. Given  integers $L>0, d\geq 1,$ let $\Lambda=[-L, -L+1, ..., 0, ...,L-1, L]^d\in \Z^d$ be a lattice cube with periodic boundary conditions. Call   $x,y\in \Lambda$  neighbours,  if $\|x-y\|=1$, where   $\|\cdot\|$  is the Euclidean norm in  $\R^d$. In this case $\nu(x)\equiv 2d$ and the  graph does not have triangles. 
\item[ b)] {\it  Mean-field model.} Given $n\geq 2$ let   $\Lambda$  be  a complete graph  with $n$ vertices. By construction, $\nu(x)\equiv n-1$ in this example and the graph {\em does have} triangles.
\end{itemize}
The following statement complements Theorem~\ref{T3} in the mean field case.
\begin{theorem}
\label{T31}
Let $\L$ be a complete graph with $n$ vertices labelled by $1, \ldots, n$, where $n\geq 1$.  If either $0<\alpha<\beta$ or  $\alpha<0<\alpha+\beta\nu$ then 
\begin{itemize}
\item[ 1)] 
${\zeta_k(t)}/{t}\to {1}/{n}$ for all $k=1,\ldots,n$ a.s.;
\item[ 2)] all components of CTMC $\xi(t)$ explode simultaneously a.s.;
\item[ 3)] a process of differences $(\zeta_1(t)-\zeta_{n}(t),\ldots, \zeta_{n-1}(t)-\zeta_{n}(t))\in\Z^{n-1}$ converges in distribution as $t\to \infty$.
\end{itemize}
\end{theorem}
Finally, Theorem~\ref{T4} below describes the long-term behaviour of the Markov chain in the case of a star graph.
\begin{theorem}\label{T4}
Given $n\geq 1$ let  $\L$ be a star graph  with $(n+1)$ vertices, i.e.\ where there is a central vertex $x$ and its neighbouring vertices $y_1, \ldots, y_n$, so that $x$ is the only neighbour for each of $y_i$, $i=1,\ldots,n$, and $x\sim y_i,\, i=1,\ldots,n$.
Then
\begin{itemize}
\item[1)] CTMC $\xi(t)$ is ergodic if and only if $\alpha<0$ and $\alpha+\beta\sqrt{n}<0$;
\item[2)] if $\alpha<0$ and $\alpha+\beta\sqrt{n}=0$ then CTMC  $\xi(t)$ is transient;
\item[3)] if  $\alpha<0$ and $\alpha+\beta\sqrt{n}>0$ then with probability $1$
$$
\frac{\zeta_{x}(t)}{t}\to \frac{n\beta+|\alpha|}{(n+1)\beta+2|\alpha|},\quad 
\frac{\zeta_{y_i}(t)}{t}\to \frac{\beta+|\alpha|}{(n+1)\beta+2|\alpha|},\ i=1,2,\dots,n,
$$
as $t\to \infty$, and hence with probability $1$ all components of CTMC $\xi(t)$ explode simultaneously;
\item[4)] if   $\alpha>0$ then  CTMC  $\xi(t)$ is explosive. Moreover, 
\begin{itemize}
\item[i)] if $\beta<\alpha$ then with probability $1$ the event~(\ref{B}) occurs  and a single component of CTMC $\xi(t)$ explodes;
\item[ii)] if  $\alpha<\beta$ then with probability $1$ the event $B_{x, y_i}$ occurs for some  $i=1,\ldots,n$, so that with probability $1$ 
only a pair of adjacent components of the CTMC explodes.
\end{itemize}
\end{itemize}
\end{theorem}

\section{Random walk in the quarter plane}\label{exp}
Let graph $\L$ be formed by two adjacent vertices. This is the simplest example of both constant degree graphs and star graphs, while the corresponding  Markov chain is equivalent to an inhomogeneous random walk on the positive quarter plane.  We will briefly comment on this particular case to illustrate some distinctive features of the model dynamics, which can be also observed in more general situations.
\begin{figure}[htbp]
\centering
\begin{tabular}{cc}
   \includegraphics[width=2.7in, height=2.2in]{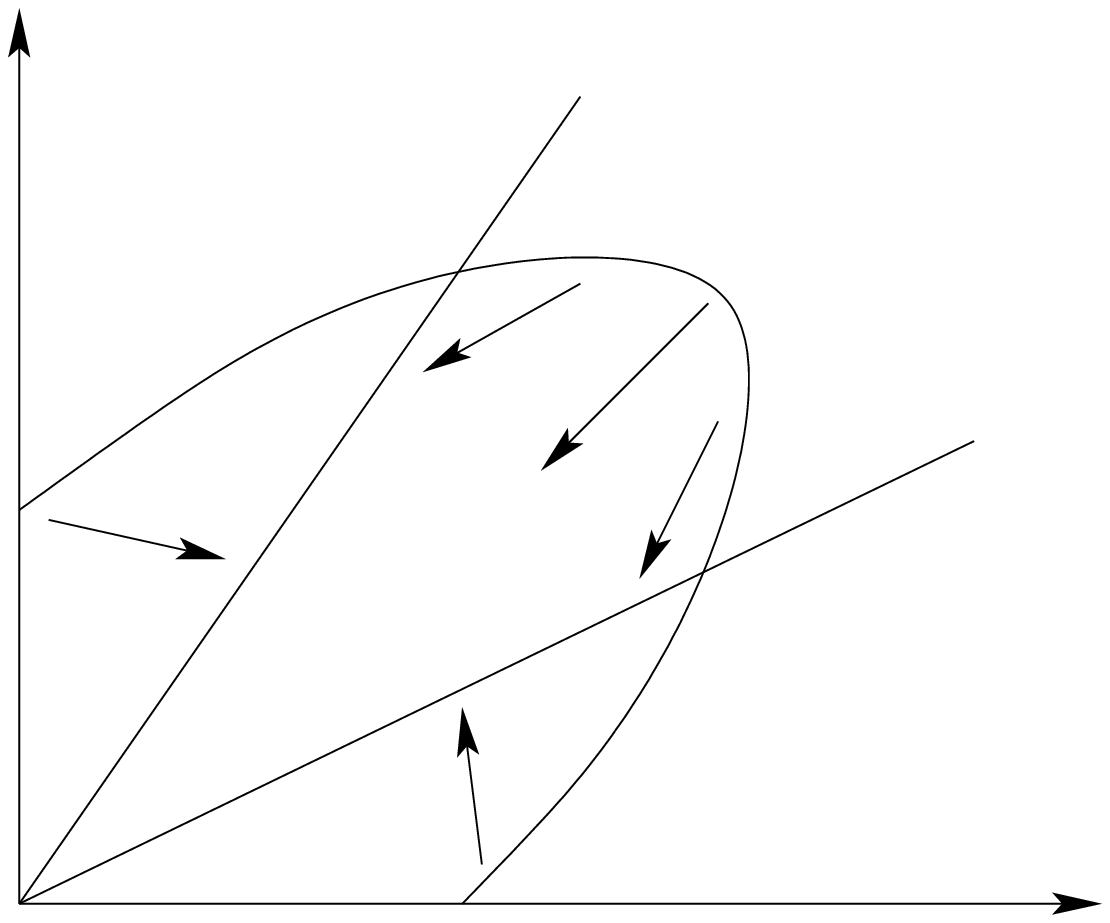}&
     \includegraphics[width=2.7in, height=2.2in]{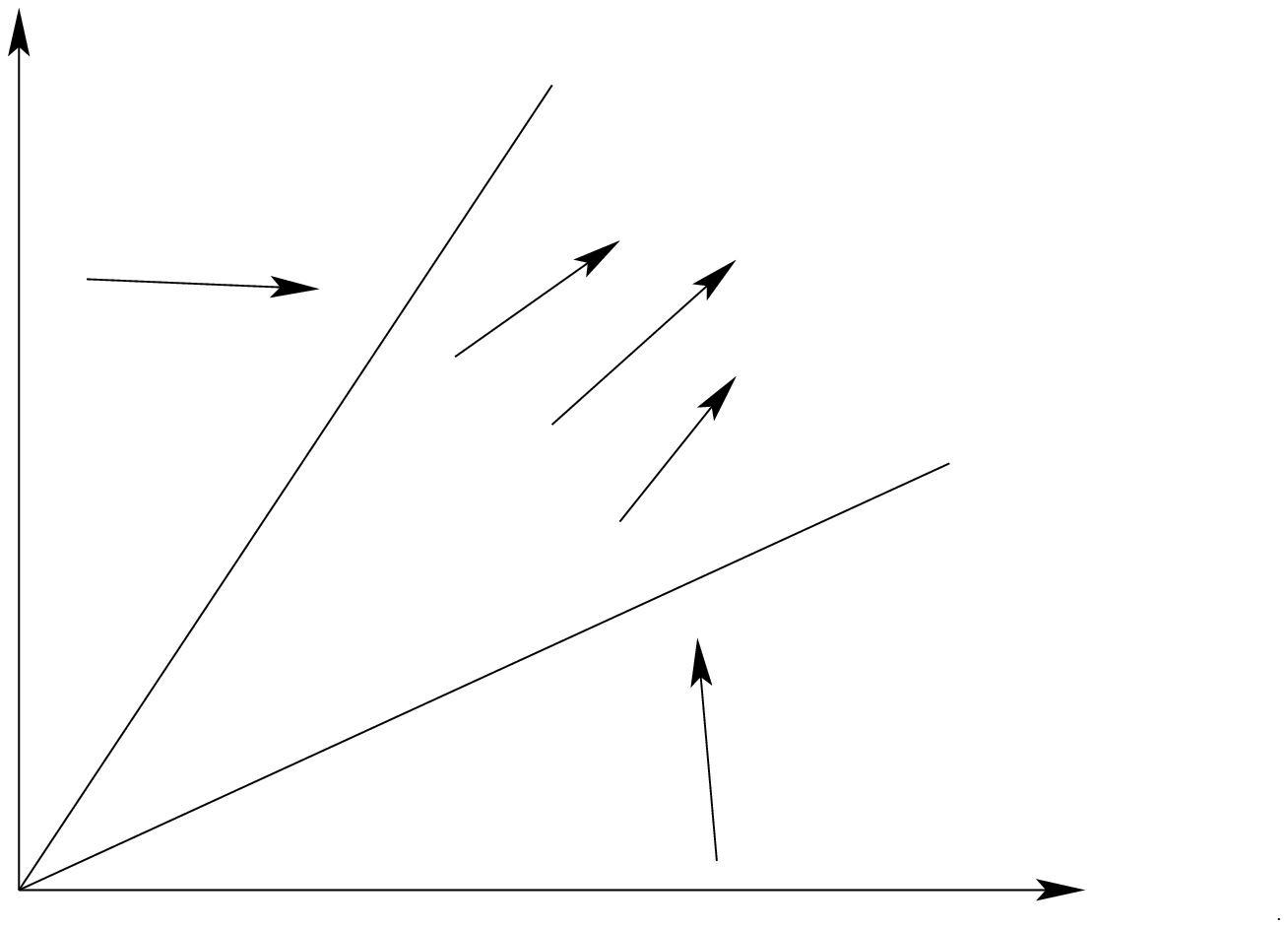} \\
\end{tabular}
\caption{{\footnotesize The vector field of mean jumps of Markov chains, $\alpha<0,\, \beta>0.$
The vertical axis is $y$ axis and the horizontal axis   $x$ axis.
the upper line is $y=-\frac{\alpha}{\beta}x$, 
the lower line  is  $y=-\frac{\beta}{\alpha}x,$ the curve is  $Q(x,y)=C$, for some $C>0$.
Right: $\alpha+\beta>0$ (transience); 
the upper line is $y=-\frac{\beta}{\alpha}x$, 
the lower line is $y=-\frac{\alpha}{\beta}x$.}}
\label{Fig1}
\end{figure}
\begin{figure}[htbp]
\centering
\begin{tabular}{cc}
   \includegraphics[width=2.7in, height=2.2in]{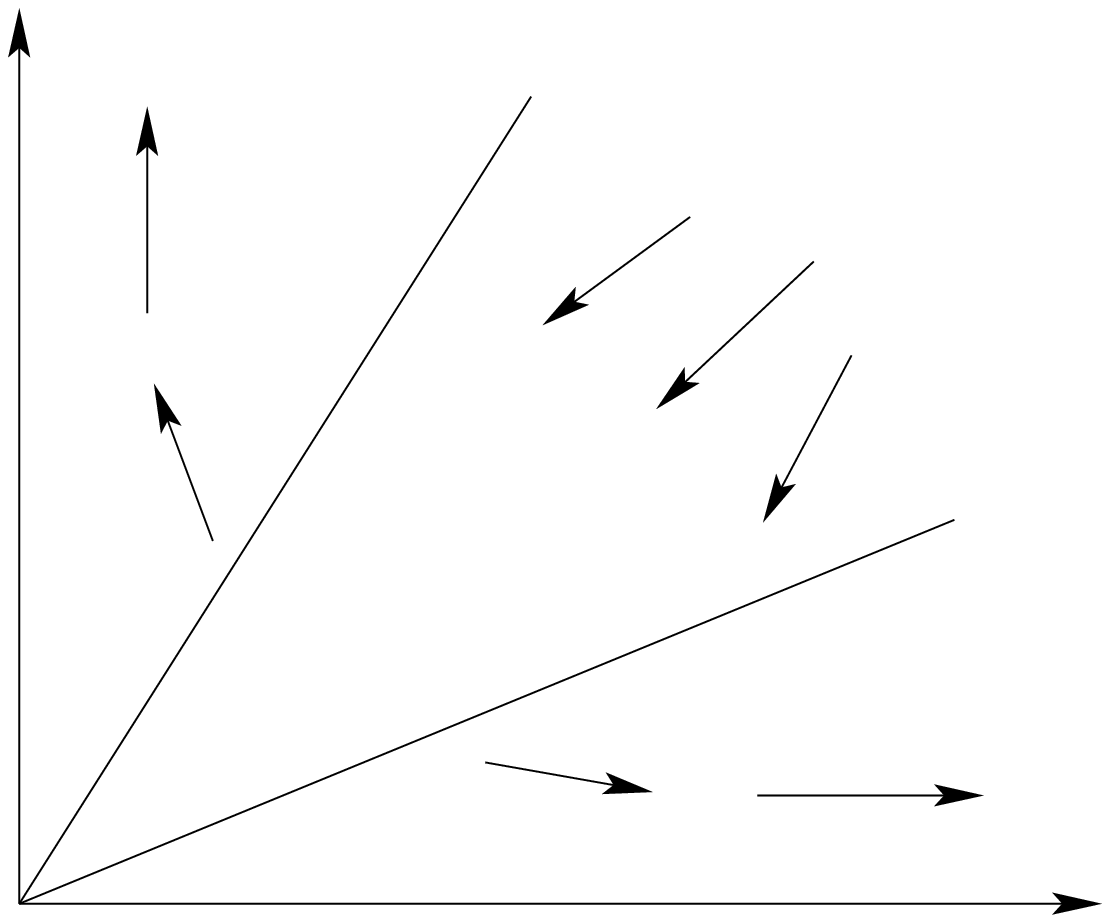}&
     \includegraphics[width=2.7in, height=2.2in]{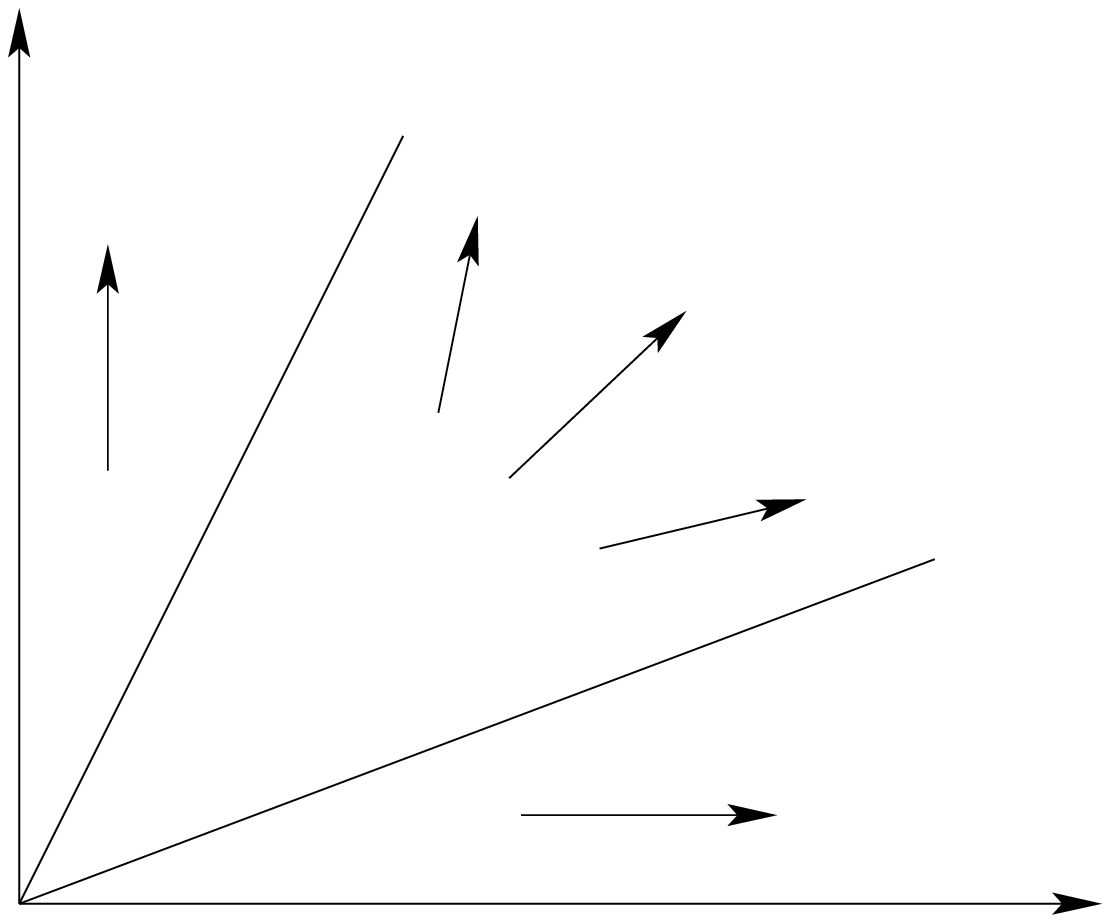} \\
\end{tabular}
\caption{{\footnotesize 
The vector field of mean jumps of Markov chains, $\alpha>0,\, \beta<0.$
 The vertical axis is  $y$ axis  and the horizontal axis is   $x$ axis.
the upper line is $y=-\frac{\beta}{\alpha}x$, 
the lower line is $y=-\frac{\alpha}{\beta}x$.
Right: 
$\alpha+\beta>0$; 
the upper line is $y=-\frac{\alpha}{\beta}x$, 
the lower line is $y=-\frac{\beta}{\alpha}x$.}}
\label{Fig2}
\end{figure}

The theorems in Section~\ref{results} imply the following results for the two-dimensional case.

1) If $\alpha<0$ and $\beta< |\alpha|$ then both CTMC $\xi(t)$ and DTMC $\zeta(t)$ are ergodic. Left part of Figure~\ref{Fig1} sketches  the vector field  of mean  jumps of the Markov chain and level curves of Lyapunov function $Q(x,y)=-\alpha x^2-\alpha y^2 - 2\beta xy$
in the ergodic case  $0<\beta<-\alpha$.

2) If $\alpha<0$ and $\alpha+\beta\geq 0$ then DTMC $\zeta(t)$ is transient, CTMC  $\xi(t)$ is explosive; moreover.
$$
\P(\zeta_1(t)=\zeta_2(t)\ \mbox{infinitely often})=1.$$
The vector field of mean jumps in the case  $\alpha<0, \alpha+\beta>0$ is illustrated by the right part of Figure~\ref{Fig1}.

3) If $\alpha>0$ then CTMC $\xi(t)$ is explosive. If, in addition, $\beta<\alpha$ then with probability $1$  a single component of DTMC will eventually grow (event~(\ref{B}) occurs). We illustrate this by the left part of Figure~\ref{Fig2} in the case $\beta<-\alpha<0$.  The right part of the same figure corresponds to the  transient/explosive  case $-\alpha<\beta<0$. If $\alpha<\beta$ then both components grow and $$
\P(\zeta_1(t)=\zeta_2(t)\,\, \mbox{infinitely often})=1.
$$

4) In the two-dimensional case we also deal with the case $\alpha=0$ and  $\beta< 0$ and show that  both CTMC  $\xi(t)$  and DTMC $\zeta(t)$ are null recurrent. Indeed, by the well-known criteria for recurrence (e.g., Theorem~2.2.1 in~\cite{FMM}) to show recurrence in both cases
 it suffices to find a positive function $f(x,y)$ such that $f(x,y)\to \infty$ as $\sqrt{x^2+y^2}\to \infty$ and for which the following inequality 
\begin{equation} \label{null}
\G f(x,y)\leq 0,
\end{equation} 
holds for all but finitely many  $(x,y)$, where $\G$ is the generator of the Markov chain, see (\ref{generator}).
 Consider a function   $f(x,y)=\log(x+y+1)$.
We will show that if the sum  $x+y$ is sufficiently large,  then the inequality~(\ref{null}) holds. Indeed, if $y=0$ then  
\begin{align*}
\G & f(x, 0)=\left(\log(x+2)-\log(x+1)\right)\left(1+e^{\beta x}\right)+\left(\log(x)-\log(x+1)\right)\\
&=\log\frac{x+2}{x+1}+\log\frac{x}{x+1}+e^{\beta x}\log\frac{x+2}{x+1}
\leq \log\left(1-\frac{1}{(x+1)^2}\right)+e^{\beta x}\leq 0,
\end{align*}
where the last inequality holds for sufficiently  large $x>0$. If both $x>0$ and $y>0$ then assuming that  $C=x+y$ is large enough, we have the following bound:
\begin{align*}
\G f(x,y)&=\left(\log(C+2)-\log(C+1)\right)\left(e^{\beta y}+e^{\beta  x}\right)
\\
&+\left(\log(C)-\log(C+1)\right)\left(1_{\{x>0\}}
+1_{\{y>0\}}\right)\\
&\leq 2\left(\log(C+2)-2\log(C+1)+\log(C)\right)=2\log\frac{C(C+2)}{(C+1)^2}\leq 0.
\end{align*}
Therefore, both Markov chains are recurrent. It is easy to see that CTMC  $\xi(t)$ cannot be positive recurrent in this case. Indeed,  had it 
 been recurrent, then its stationary distribution would be given by formula~(\ref{inv}), but the latter is impossible, since 
$$
\sum\limits_{x,y\in Z_{+}} e^{\beta  xy}=\infty
$$ 
for all $\beta$. This also yields that DTMC $\zeta(t)$ cannot be positive recurrent as well. 

\section{Proofs}
\label{proofs}

\subsection{Proof of Theorem~\ref{T1}}

{\it Proof of  Part 1) of Theorem~\ref{T1}.}  Notice first that if~$\alpha<0$ and~$\beta\leq 0$ then
   the stationary distribution~(\ref{inv}) is  well-defined and  the Markov chain $\xi(t)$ is ergodic.

 We will show now that CTMC does not explode, if $\alpha<0,\, \beta>0$ and $\alpha+\beta\max_{x\in \L}\nu(x)\leq 0$. 
Define
$$
\tau_N=\min\left\{t: \max\limits_{x\in \L}\xi_x(t)=N\right\}.
$$
 It is obvious that the Markov chain is explosive if and only if 
$$
\P\left(\lim\limits_{N\to\infty} \tau_N<\infty\right)>0,
$$
but the latter  cannot happen. 
Indeed, given $\xi(t)=\xi$ let $x\in \L$ be such that 
 $\xi_x=\max_{y\in \L} \xi_y$. Then 
$$
U(x, \xi)=\alpha \xi_x+\beta\phi(x, \xi)\leq (\alpha +\beta \nu(x) )\xi_x\leq \left(\alpha+\beta\max_{x\in \L}\nu(x)\right)\xi_x\leq 0.
$$
Therefore the waiting times $\tau_{N+1}-\tau_N$ are stochastically larger than  exponentially distributed independent random variables with  parameter $(2|\L|)^{-1}$; as a result,  the limit $\lim\limits_{N\to\infty} \tau_N$ is infinite with probability $1$ and thus the chain does not explode.

Let us finally show that if 
\begin{equation}
\label{strict}
\alpha<0,\, \beta>0, \alpha+\beta\max_{x\in \L}\nu(x)< 0,
\end{equation}
 then $Z_{\alpha, \beta, \L}<\infty$ 
and consequently the stationary probability distribution is well defined. 
It is easy to see that $Z_{\alpha, \beta, \L}<\infty$ if and only if
\begin{equation}
\label{summ}
\sum\limits_{\xi\in \Z_{+}^{\L}}\exp(-Q(\xi)/2)<\infty.
\end{equation}
Consider a symmetric matrix   $A_Q=(a_{xy})_{x,y\in \L}$   determining  the quadratic form $Q$, i.e.
\begin{equation}
\label{AQ}
Q(u)=(A_Qu, u),\quad u\in \R^{|\L|}.
\end{equation}
It is easy to see that 
$a_{xx}=-\alpha,\, a_{xy}=-\beta,$ if $y\sim x$ and $a_{xy}=0$ otherwise. Therefore for all $x\in \L$
$$|a_{xx}|-\sum\limits_{y\neq x} |a_{xy}|=-\alpha-\beta\nu(x)\geq -\alpha-\beta\max_{x\in \L}\nu(x)>0,$$
 because of (\ref{strict}).
In other words, matrix $A_Q$ is strictly  diagonally dominant with positive diagonal  entries and hence, by  standard algebra, 
 is positive definite.  One can also observe positive definiteness of $A_Q$ in the  case under consideration from representation~(\ref{Qform1}). 
 Positive definiteness of $A_Q$ implies that 
$$
\int\limits_{\R^{|\L|}}e^{-(A_Qu,u)/2}du=\frac{(2\pi)^{\frac{|\L|}{2}}}{\sqrt{\mbox{det}(A_Q)}}<\infty,
$$
which, in turn, implies~(\ref{summ}), so the  stationary probability distribution is well defined as claimed.

{\it Proof of Part 2) of Theorem~\ref{T1}.}
The Markov chain $\xi(t)$ cannot be ergodic if $\alpha\geq 0$. Indeed, fix $x\in \L$ and  define  a set of configurations
$D_x=\{\xi: \xi_x\geq 0,\, \xi_y=0,\, y\neq x\}$.
It is easy to see that 
$$Z_{\alpha, \beta, \L}\geq \sum\limits_{\xi\in D_x} e^{W(\xi)}= \sum\limits_{k=0}^{\infty}e^{\alpha k(k-1)/2}=\infty,$$
and  the stationary distribution does not exist.  

\paragraph{Function $Q$ as the Lyapunov function for Foster criteria.}
Observe that ergodicity of the Markov chain in Part~1)  of Theorem~\ref{T1} can be shown by using Foster criteria for ergodicity of a countable Markov chain. We skip the easy case, when both  $\alpha<0$ and $\beta<0$ and show that if $\alpha<0, \beta>0$ and $\alpha+\beta\max_{x\in \L}\nu(x)<0$  the function $Q$  serves as the corresponding Lyapunov function. Indeed,  the equation~(\ref{Qform1}) yields that  $Q(\xi)>0$ for all $\xi\in \Z_{+}^{\L}$ outside the origin (i.e., $\xi\neq 0$) and that $Q(\xi)\to \infty$ as $\sum_{x\in\L}\xi_x^2\to \infty$.  Recall that  $\G$ is the generator  (defined by (\ref{generator}))
of the Markov chain.
 We fix some $\eps>0$ and show that 
\begin{equation}
\label{fost1}
\G Q(\xi)\leq -\eps,
\end{equation}
provided that $S(\xi)=\sum_{x\in \L}\xi_x\geq C$, where $C=C(\eps)$ is sufficiently large. It is easy to see that 
\begin{equation}
\label{fost11}
\G Q(\xi)=\sum\limits_{x\in \L}(-\alpha-2U(x, \xi))e^{U(x, \xi)}+\sum\limits_{x\in\L}(-\alpha+2U(x, \xi))1_{\{\xi_x>0\}},
\end{equation}
where $U(x, \xi)$ is defined by equation~(\ref{U}). Sums  in~(\ref{fost11}) can be respectively  bounded   as follows
$$
\sum\limits_{x\in \L}(-\alpha-2U(x, \xi))e^{U(x, \xi)}\leq 
|\L|\max_{u\in \R}(-\alpha+2u)e^{-u}=2|\L|e^{\frac{-\alpha-2}{2}},
$$
and   
\begin{align*}
\sum\limits_{x\in \L}(-\alpha+2U(x, \xi))1_{\{\xi_x>0\}}&\leq 
\sum\limits_{x\in \L}(-\alpha+2U(x, \xi))
=-\alpha|\L| +2\sum\limits_{x\in \L}(\alpha+\beta\nu(x))\xi_x\\
&\leq -\alpha|\L| +2(\alpha+\beta\max_{x\in \L}\nu(x))S(\xi)
\\
&\leq -\alpha|\L|+2C(\alpha+\beta\max_{x\in \L}\nu(x)),
\end{align*}
where we used  the equation~(\ref{Usum}) to get the second display. Thus the LHS of~(\ref{fost1}) is bounded by the following quantity 
$$
2|\L|e^{\frac{-\alpha-2}{2}}-\alpha|\L|+2C(\alpha+\beta\max_{x\in \L}\nu(x)),
$$
which  is less than $-\eps$,  if $C>0$ is sufficiently large. The inequality~(\ref{fost1}) allows to apply  Foster criteria of ergodicity (Theorem~2.2.3 in~\cite{FMM}) of a countable  Markov chain.

\subsection{Proof of Theorem~\ref{T2}} 

We start with showing that  there exists a $\delta'>0$ such that
\begin{equation}
\label{P01}
\P\left(B|\zeta(t)=\zeta\right)>\delta',
\end{equation}
for all  $t$ and $\zeta$. 
Given $\zeta\in \Z_{+}^{\L}$ define 
\begin{align*}
M_{\zeta}=\max_{x\in \L}U(x, \zeta)
\ \mbox{ and } \
D_{\zeta}=\left\{x\in \L: U(x, \zeta)=M_{\zeta}\right\}.
\end{align*}
Given $\{\alpha, \beta\}$ there can be two different cases. 
\begin{enumerate}
\item[1)] A finite connected graph $\L$ is such that 
\begin{equation}
\label{U>01}
M_{\zeta}\geq 0\,\, \mbox{for all}\,\, \zeta\in \Z_{+}^{\L}.
\end{equation}
We say in this case that $\L$ is a type I graph. 
\item[2)] The set of configurations
\begin{equation}
\label{K}
{\cal K}=\{\zeta: M_{\zeta}<0\},
\end{equation}
is not empty, then    we say that  $\L$ is a type II graph.
\end{enumerate}
Let us  consider some examples before proceeding further.
It is obvious that if both  $\alpha$ and $\beta$ are positive, then any graph is a type I graph. 
Also, if $\alpha>0>\beta$ and  $\alpha+\beta\max_{x\in \L}\nu(x)\geq 0$ then  
 for every  $x\in \L$ such that $\zeta_x=N=\max_{y\in \L}\zeta_y$ the following inequality holds
$$
U(x, \zeta)=\alpha N+\beta\phi(x, \zeta) \geq N\left(\alpha+\beta\max_{x\in \L}\nu(x)\right)\geq 0,
$$
hence, $\L$ is a  type I graph.

Consider also two main examples  of type II graphs. The first one is when  $\alpha>0>\alpha+\beta\nu$, and  $\L$ is a constant vertex degree graph with $\nu(x)\equiv \nu$.
In this case ${\cal K}$  is a non-empty set of configurations that belongs to   intersection of hyperplanes $\{\zeta\in \R_{+}^{\L}: \alpha \zeta_x +\beta\phi(x, \zeta)<0,\, x \in \L\}$. The  second example is when $\alpha>0>\alpha+\beta\sqrt{\nu}$ and $\L$ is a star graph with $n+1$ vertices. In this case the set  ${\cal K}$ is the subset of intersection of hyperplanes  $\{\zeta\in \R_{+}^{\L}: \alpha \zeta_x +\beta\phi(x, \zeta)<0,\, x \in \L\}$. Notice also that  if  set ${\cal K}$  not empty, then it is  infinite.  Indeed, if $\zeta\in {\cal K}$ then $a\zeta\in {\cal K}$ for any $a\in \Z_{+}$.

Let us continue with the proof of~(\ref{P01}).
 First, given  a pair  $\{\alpha, \beta\}$  suppose that $\L$ is a type I graph. For a given $x\in \L$ define the following event
$$
B_{t, x}=\{\zeta_{x}(s+1)=\zeta_{x}(s)+1,\, \zeta_y(s)=\zeta_y(t),\,\, \mbox{for}\,\,  y\neq x\,\, \mbox{and}\,\, s\geq t\}.
$$
Trivially, $B_{t,x}\subset B$. We are going to show that for any~$\zeta$  and~$x\in D_{\zeta}$
$$
\P(B_{t, x}\| \zeta(t)=\zeta)>\delta'>0, 
$$
where $\delta'$  might depend only on parameters $\alpha, \beta$ and graph $\L$. Given  $x\in \L$ and $\zeta\in \Z_{+}^{\L}$  denote 
$$
R(x, \zeta)=\sum\limits_{y\in \L} e^{U(y, \zeta)}-\left(e^{U(x, \zeta)}+\sum\limits_{y\sim x}e^{U(y, \zeta)}\right).
$$
If $x\in D_{\zeta}$ then  
\begin{equation}
\label{R1}
R(x, \zeta)e^{-U(x, \zeta)}=R(x, \zeta)e^{-M_{\zeta}}\leq (|\L|-\nu(x)-1)<|\L|,
\end{equation}
for all $\zeta\in \Z_{+}^{\L}$. Given~$x\in D_{\zeta}$  we have that    
\begin{align*}
\P\left( B_{t, x}| \zeta\right)&= \prod\limits_{k=0}^{\infty}
\frac{e^{M_{\zeta}+\alpha k}}{e^{M_{\zeta}+\alpha k}+{\displaystyle\sum_{y\sim x}}e^{U(y, \zeta) +\beta k}+R(x, \zeta) +\sum_{y\in \L}1_{\{\zeta_y>0\}}}\\
&= \prod\limits_{k=0}^{\infty}
\frac{1}{1+{\displaystyle\sum_{y\sim x}}e^{U(y, \zeta)-M_{\zeta} -(\alpha-\beta) k}+\left[R(x, \zeta) 
 +{\displaystyle\sum_{y\in \L}}1_{\{\zeta_y>0\}}\right]e^{-M_{\zeta}-\alpha k}},
\end{align*}
for all $\zeta\in \Z_{+}^{\L}$. It is easy to see that by  choice of $x$  we have
$$
\sum_{y\sim x}e^{U(y, \zeta)-M_{\zeta} -(\alpha-\beta) k}\leq  e^{-(\alpha-\beta)k}\max_{y\in \L}\nu(y).
$$
Also, using~(\ref{U>01}) and~(\ref{R1}) we get that 
\begin{equation}
\label{R2}
\left(R(x, \zeta) +\sum_{y\in \L}1_{\{\zeta_y>0\}}\right)e^{-M_{\zeta}-\alpha k}\leq 2|\L|e^{-\alpha k}.
\end{equation}
Therefore, we obtain the following bound
\begin{equation}\label{delta1}
\P\left( B_{t, x}\| \zeta\right)\geq \prod\limits_{k=0}^{\infty}
\frac{1}{1+e^{-(\alpha-\beta)k}\max_{y\in \L}\nu(y)+2|\L|e^{-\alpha k}}=\delta'>0.
\end{equation} 
The preceding display implies  bound~(\ref{P01}) in the case of a type I graph. 

Suppose now that  given a pair $\alpha, \beta$  satisfying   the conditions of the theorem, $\L$ is a type II graph.
Fix some $\eps>0$ and suppose that
$\zeta\in {\cal K}_{\eps}=\{\zeta:  M_{\zeta}\geq -\eps\}$.
Given  $x\in D_{\zeta}$  one can repeat, with a minor change,  the same argument
which led to  bound~(\ref{delta1}). The only difference now  is that  the inequality $M_{\zeta}\geq -\eps$ yields 
constant $(1+e^{\eps})|\L|e^{-\alpha k}$ in the right  side of (\ref{R2})  (instead of  $2|\L|e^{-\alpha k}$)
and it results in  a different $\delta''\neq \delta'$  such that 
$$
\P(B_{t,x}\| \zeta(t)=\zeta)>\delta''>0.
$$
Consider  the case, when $\zeta\in {\cal K}^c_{\eps}=\{\zeta: M_{\zeta}<-\eps\}$.
Define a stopping time
$$
\tau=\min\{t: \zeta(t)\in {\cal K}_{\eps}\}.
$$
It is easy to see that  $\P(\tau<\infty|\zeta)=1$ for all $\zeta\in {\cal K}^c_{\eps}$. 
 Indeed, define $F(\zeta)=\|\zeta\|^2$, where $\|\zeta\|$ is Euclidean norm in $\R^{\L}$.
 A direct computation gives that there exists  some $\eps'>0$ such that 
$$\E(F(\zeta(t+1))-F(\zeta(t))|\zeta(t)=\zeta)\leq -\eps'$$
for all $\zeta\in {\cal K}^c_{\eps}$.  The inequality in the preceding display yields, by a  standard argument (see, for example, Theorem~2.1.1 in~\cite{FMM}), that starting from inside of  ${\cal K}^c_{\eps}$ the process  $\zeta(t)$  reaches the set $\{\zeta: \|\zeta\|\leq C\}\cap{\cal K}^c_{\eps}$ in a finite  time unless it first  exits  the set ${\cal K}^c_{\eps}$. Suppose the process reaches set $\{\zeta: \|\zeta\|\leq C\}$ 
before exiting  ${\cal K}_{\eps}^c$. 
 DTMC $\zeta(t)$  transits to ${\cal K}_{\eps}$ from  $\{\zeta: \|\zeta\|\leq C\}$ with positive probability
due to irreducibility of DTMC and finiteness of set $\{\zeta: \|\zeta\|\leq C\}$. 
This implies that  DTMC $\zeta(t)$ reaches  set ${\cal K}_{\eps}$ with probability $1$ and we can apply the same argument 
as in the case of type I graph to obtain  the following  bound
$$
\P(B\|\zeta(t)=\zeta)=\P\left(B, \tau<\infty\|\zeta\right)\geq \min_{\zeta\in  {\cal K}^c, x\in D_{\zeta} }\P\left(B_{\tau, x}\|\zeta\right)>\delta''>0.
$$ 
for all $\zeta\in {\cal K}^c_{\eps}$.
This completes the proof of bound~(\ref{P01}).

Bound~(\ref{P01})  implies  that 
\begin{equation}
\label{>01}
\E\left(1_{B}\|{\cal F}_t\right)>\delta',
\end{equation}
where ${\cal F}_t=\sigma\{\zeta_0, \ldots, \zeta_t\}$ is the $\sigma-$algebra of events generated by DTMC  up to time moment $t$. 
Since $B\in {\cal F}_{\infty}=\sigma\{{\cal F}_t,\, t\geq 0\}$ we get by Levi $0-1$ law that 
$$
\E(1_{B}\|{\cal F}_t)\to \E(1_{B}\|{\cal F}_{\infty})=1_{B},\,\, \mbox{as}\,\, t\to \infty.
$$
By~(\ref{>01}) the right hand side  of the preceding display is positive. Therefore, it must be equal  to $1$, hence,  $\P(B)=1$.

Thus, eventually only a single component of the embedded chain continues to evolve by jumping up without jumping down. In the continuous time setting the only  growing component  evolves  eventually as a pure birth process with exponentially growing birth rates. The latter  process is explosive and, hence,  CTMC  $\xi(t)$ is explosive, where with probability $1$ only a single component explodes.

Notice also that under the assumptions of the theorem  with probability one a typical  trajectory of DTMC $\zeta(t)$   returns 
to set    ${\cal K}^c_{\eps}$ only a finite number of times in the case of type II graph.

\subsection{Proof of Theorem~\ref{T.no.triangle}}
We start with  the following lemma.
\begin{lemma}
\label{L1}
Let $0<\alpha<\beta$. Suppose that $x_1$ and $x_2$ are two vertices of $\Lambda$ such that (1) $x_1\sim x_2$; (2) there is no $y$ such that $y\sim x_1$ and $y\sim x_2$ at the same time; (3) at some time $s$ the configuration of the DTMC is such that $u_1=U(x_1,\zeta(s))$ is the largest potential on the whole graph and $u_2=U(x_2,\zeta(s))$  is the largest potential among all the neighbours of $x_1$. Then, with a positive probability depending on $\alpha$, $\beta$ and $\Lambda$ only, the following events simultaneously occur
\bn
\zeta_y(t)&=\zeta_y(s)\text{ for all }y\notin\{x_1,x_2\}
\text{ and all }t=s,s+1,s+2,\dots; \nonumber\\
\lim_{t\to\infty} \frac{\zeta_{x_1}(t)}{t}&
= \lim_{t\to\infty} \frac{\zeta_{x_2}(t)}{t}=\frac{1}{2}.
\en
\end{lemma}
{\it Proof of Lemma~\ref{L1}.}
Observe that every time when the component at $x_1$ increases by $1$, the potential at $x_1$ increases by $\alpha$ while the potential at each of the neighbours of $x_1$ increases by $\beta$, therefore the potential at $x_2$ remains the largest among the neighbours of $x_1$. At the same time the difference between the potentials at $x_1$ and $x_2$ decreases by $\delta:=\beta-\alpha>0$.

Let $k=\left\lfloor\frac{u_1-u_2}{\delta}\right\rfloor$ where $\lfloor a\rfloor$ denotes the integer part of $a\in\R$. W.l.o.g.\  assume that $k$ is even. Denoting by $\nu_1=\nu(x_1)$ the degree of vertex $x_1$ (and $\nu_2=\nu(x_2)$ respectively), we obtain that the probability that during the times $t=s,s+1,...s+k$ only the component at $x_1$ increases is larger than
\bn
p_1&=\prod_{i=0}^{k}
\frac{e^{u_1+i\alpha}}{e^{u_1+i\alpha}+\nu_1 e^{u_2+i\beta}+(|\Lambda|-\nu_1)e^{u_1}}
\\
&=\prod_{i=0}^{k}
 \frac{1}{1+\nu_1 e^{-[u_1-u_2]+i(\beta-\alpha)}+(|  \Lambda|-\nu_1)e^{-i\alpha}}
\ge 
\prod_{i=0}^{k} \frac{1}
{1+\nu_1 e^{-(k-i)\delta}+|\Lambda|e^{-i\alpha}}
\\ &  \ge 
\prod_{i=0}^{k/2} \frac{1}
{1+\nu_1 e^{-k\delta/2}+|\Lambda|e^{-i\alpha}} \times
\prod_{j=0}^{k/2} \frac{1}
{1+\nu_1 e^{-j\delta}+|\Lambda|e^{-k\alpha/2}}
\\
&
\ge 
\left(\prod_{i=0}^{k/2} \frac{1}
{1+(\nu_1+|\L|)( e^{-i\delta}+e^{-i\alpha})}\right)^2
=C_1(|\Lambda|,\alpha,\beta)>0.
\en
Consequently, by time $s+k$ we have have $-\delta<U(x_2,\zeta(s+k))-U(x_1,\zeta(s+k))\le 0$ with probability at least $p_1$.

From now on assume w.l.o.g.\  that actually $u_2\in (u_1-\delta,u_1]$ already at time $s$. Let $m_i(t)$, $i=1,2$ be the number of times $x_i$ was chosen during the times $s+1,s+2,\dots,s+t$. Define the events
\bn
A_k'&=\{\zeta_y(s+i)=\zeta_y(s)\text{ for all }y\notin\{x_1,x_2\},\ i=1,2,\dots,2k^2 \}
\\
A_k''&=\left\{|m_1(2k^2)-m_2(2k^2)|\le 2k\right\}
\\
A_k&=A_k'\cap A_k''.
\en
Then under $A_k$ we have $m_1(2k^2)+m_2(2k^2)=2k^2$ and $|m_i(2k^2)-k^2|\le k$ for $i=1,2$, so $\P(A_{k+1}' | A_k)$ is no less than
\bn
& 
\prod_{i=0}^{4k+1}
\frac{e^{U(x_1,\zeta(s+2k^2+i))}+e^{U(x_2,\zeta(s+2k^2+i))}}
{e^{U(x_1,\zeta(s+2k^2+i))}+e^{U(x_2,\zeta(s+2k^2+i))}
+[\nu_1+\nu_2]e^{u_2+\beta (k^2+6k)}+ [|\Lambda|-\nu_1-\nu_2]e^{u_2}}
\\
&\ge
\prod_{i=0}^{4k+1} \frac 1{1+|\Lambda|e^{(7\beta+\alpha)k-\alpha k^2}}
\ge 1-C_2(|\Lambda|,\alpha) e^{-k}
\en
since $U(x_1,\zeta(s+2k^2+i))\ge u_1+\alpha(k^2-k)+\beta(k^2-k)$, and the potential at any $y\sim x_1$ or $\sim x_2$ is bounded by $u_2+\beta(k^2+k+(4k+1))\le u_2+\beta(k^2+6k)$. To estimate
$\P(A_{k+1}'' \| A_k)$ observe that whenever $m_1(j)>m_2(j)+1$ the potential at $x_2$ is larger, and the similar statement holds if one swaps $1$ and $2$. Now, there are two possibilities at time $j=s+2k^2$: (a)~$|m_1(2k^2)-m_2(2k^2)|\le 1.5 k$ and (b)~$|m_1(2k^2)-m_2(2k^2)|> 1.5k$. 

In case~(a), the difference $|m_1(j)-m_2(j)|$ can be majorized by the distance to the origin of the simple symmetric random walk on $\Z^1$. In particular, the probability that during $4k+2$ steps it is further than $k^{2/3}$ from the starting point is bounded by $c_3 e^{-k^{1/6}}$ where $c_3$ is some constant. As a result,, with probability at least $1-c_3 e^{-k^{1/6}}$ we have $|m_1(2(k+1)^2)-m_2(2(k+1)^2)|< 1.5k+k^{2/3}<2(k+1)$ and $A''_{k+1}$ occurs. 

On the other hand, in case~(b) we have $1.5k<|m_1(2k^2)-m_2(2k^2)|\le 2k$, hence the potential at the larger $x_i$  in the pair $\{x_1,x_2\}$ is much smaller than the potential at the smaller $x$ in this pair. Consequently, for the next $k$ steps the probability to increase the larger component, divided by the probability to increase the smaller component, is bounded above by $e^{-\delta k/2}$, so we can couple $|m_1(j)-m_2(j)|$  with an asymmetric simple random walk on $\Z^1$ with the drift towards the origin. As a result, we obtain that with probability at least $1-e^{-c_4 k}$ during the times $t=s+2k^2+i$, $i=1,\dots,k$, the distance between $m_1$ and $m_2$ decreases at least by $k/2$,  bringing it to the value less than $2k-(k/2)=1.5k$, and thus to case~(a). Therefore,
\bn
\P(A_{k+1}'' \| A_k)&\ge 1-C_3 e^{-k^{1/6}}-e^{-C_4 k}.
\en
Combining the above inequalities yields
\bnn\label{eqpaa}
\P(A_{k+1} \| A_k)&\ge 1-C_3 e^{-k^{1/6}}-e^{-C_4 k}-C_2(|\Lambda|,\alpha) e^{-k}.
\enn
Since the product of the terms on the RHS of~(\ref{eqpaa}) over all large enough $k$ is positive, the statement of the lemma follows.

\vskip 3mm

Now note that at any moment of time $s$ there is a vertex $x_1$ with the largest potential. Because of our assumption it satisfies the conditions of Lemma~\ref{L1} for {\it some} neighbour $x_2$. Hence, Theorem~\ref{T.no.triangle} follows from the Levy 0--1 law.

\subsection{ Proof of Theorem~\ref{T3}}
{\it Proof of Part 1) of Theorem~\ref{T3}.}
Non-ergodicity in the case $\alpha\geq 0$ and ergodicity in the case $\alpha<0,\, \alpha+\beta\nu<0$ are implied by Theorem~\ref{T1}.
 If  $\alpha<0,\, \alpha+\beta\nu\geq 0$ then, using equations~(\ref{QS}) and~(\ref{Qform1}),  we get the  the following bound
$$
Z_{\alpha, \beta, \L}\geq \sum\limits_{\xi\in Z_{+}^{\L}} e^{W(\xi)}1_{\{\xi: \xi_x=\xi_y,\, \forall\,  x, y\in  \L\}}=
\sum\limits_{k=1}^{\infty}e^{|\L|((\alpha+\beta\nu)k^2-\alpha k)/2}=\infty,$$ 
which means that the stationary distribution does not exist in this case and, hence, the  CMTC is not ergodic.

Notice, in addition,  that if $\L$ is a  constant vertex degree graph then $(-\alpha-\beta\nu)$  is the eigenvalue of~$A_Q$  with the corresponding eigenvector  $(1,\ldots,1)$ and, hence, the function $\exp(-Q(\xi)/2)$ is not summable in the direction of this eigenvector, provided that $-\alpha-\beta\nu\geq 0$. Furthermore,  if  $\alpha<0, \beta>0$ then $-\alpha-\beta\nu$ is the minimal eigenvalue of $A_Q$, since all eigenvalues of matrix $A_Q$ lie,  by Gershgorin circle theorem, within the closed interval $[-\alpha-|\beta|\nu, -\alpha+|\beta|\nu]$.  

Further, it is rather straightforward to compute  the characteristic polynomial of matrix  $A_Q$ in the case of the mean-field model with $n$ vertices (complete graph with $n$ vertices). This polynomial  is 
$$
(-1)^{n-1}(\alpha-\beta+\mu)^{n-1}(-\alpha-(n-1)\beta-\mu),
$$  
and analysis of the eigenvalues yields the same results for a complete graph.
\paragraph{Proofs for non-ergodic cases.}
\begin{lemma}
\label{L3}
Let $\L$ be a finite connected graph with the constant vertex degree  $\nu(x)\equiv \nu$. If $\alpha+\beta\nu>0$ then  with probability $1$ there exists 
a time moment $\tau<\infty$ such that for all $t\geq \tau$  none of the components of DTMC $\zeta(t)$  decreases. 
\end{lemma}
{\it Proof of Lemma~\ref{L3}.}
Recall that  $U(x, \zeta)$ is the quantity  defined by equation~(\ref{U}) and  the quantity $S(\zeta)$ is defined by~(\ref{S}).
Since $\alpha+\beta\nu>0$  equation~(\ref{Usum}) implies  that for all $\zeta$
\begin{equation}
\label{MaxU}
 \max_{x\in \L} U(x,\zeta)  \ge C_1 S(\zeta),
\end{equation}
where $C_1=(\alpha+\beta\nu)/|\L|$. Using this  bound for the maximal potential we get  the following inequality 
\begin{align*}
\P(S(\zeta(t+1))=S(\zeta(t))+1\|\zeta(t)=\zeta)&=
1-\frac{\sum_{x\in \L} 1_{\{\zeta_x>0\}}}{\sum_{x\in \L}e^{U(x,\zeta)}+\sum_{x\in \L}1_{\{\zeta_x>0\}}}\\
&\geq
1-\frac{|\L|}{{\displaystyle\max_{x\in \L}}\ e^{U(x,\zeta)}}\geq 1-|\L|e^{-C_1 S(\zeta)}.
\end{align*}
Therefore, if $
D_s=\{\mbox{none of the components ever decreases after time }s\}$, then 
\begin{align}
\label{eqDs}
\P(D_s\|\zeta(s)=\zeta)\geq \prod\limits_{t=s}^{\infty}\left(1- C_2(\zeta)e^{-C_1 (t-s)}\right)
=1-o(S(\zeta))
\end{align}
where $C_2(\zeta)=|\L|e^{-C_1 S(\zeta)}$ and $o(S(\zeta))\to 0$ as $S(\zeta)\to \infty$.
 Since for any $N>0$ the set of configurations $\{\zeta:\ S(\zeta)\geq N\}$ is finite and the Markov chain is irreducible, for each $N=1,2,\dots,$  we can define
$\tau_N=\min\{t:\ S(\zeta(t))= N\}<\infty$. As
$\P(D_{\tau_N})\to 1$, by continuity of probability $\P(\cup_{N} D_{\tau_N})=1$ and hence there exists $N$ such that after time $\tau=\tau_N$ the only changes in the system are increases of the components. This finishes the proof of  Lemma \ref{L3}.

\bigskip 

It is obvious that  Lemma~\ref{L3} implies  transience of the  DTMC in the case $\alpha+\beta\nu>0$.
Nevertheless we would like to provide another lemma (Lemma~\ref{L2} below) that ensures   transience  in this case. This lemma 
takes into account   the geometry of  mean jumps and formalizes intuition which can be inferred  from, for example, 
 right images in Figures \ref{Fig1} and \ref{Fig2}.
Besides, it  provides an idea for proving transience in the case $\alpha+\beta\nu=0$ (see Lemma~\ref{L4} below).  
\begin{lemma}
\label{L2}
Let $\L$ be a finite connected graph with the constant vertex degree  $\nu(x)\equiv \nu$.
If $\alpha+\beta\nu>0$, then  for any $0<\eps<1$ the following bound holds
 \begin{equation}
\label{eps}
\E(S(\zeta(t+1))-S(\zeta(t))\|\zeta(t)=\zeta)\geq \eps,
\end{equation}
provided that  $S(\zeta)>C_1=\frac{2|\L|\eps}{(1-\eps)(\alpha+\beta\nu)}$.
\end{lemma}
{\it Proof of Lemma~\ref{L2}.} 
It is easy to see  that inequality (\ref{eps}) is equivalent to the following one 
\begin{equation}
\label{J}
 {\rm J}(\zeta, \eps):=\sum\limits_{x\in \L}\left(\delta(\eps) e^{U(x, \zeta)}-1_{\{\zeta_x>0\}}\right) \geq 0,
\end{equation}
where 
\begin{equation}
\label{delta}
\delta(\eps)=\frac{1-\eps}{1+\eps}.
\end{equation}
 Using subsequently  inequality  $1_{\{\zeta_x>0\}}\leq 1$, equation~(\ref{Usum})  and inequality $e^u>1+u,\, u\in \R$,
we obtain 
\begin{align*}
{\rm J}(\zeta, \eps)&\geq \sum\limits_{x\in \L}\left(\delta(\eps) e^{U(x, \zeta)}-1\right)\\
&\geq  \sum\limits_{x\in \L}\left(\delta(\eps) (1+U(x, \zeta))-1\right)\\
&= \delta(\eps)(\alpha+\beta\nu)S(\zeta)-(1+\delta(\eps))|\L|>0, 
\end{align*}
provided that $S(\zeta)>C_1=\frac{(1+\delta(\eps))|\L|}{\delta(\eps)(\alpha+\beta\nu)}$. 
Notice that it is also possible to  use inequality between the arithmetical and geometric means and  equation~(\ref{Usum})
in order  to obtain that 
$$\sum\limits_{x\in \L}e^{U(x, \zeta)}\geq |\L|e^{\frac{(\alpha+\beta\nu)S(\zeta)}{|\L|}}$$
and to arrive to a similar result (provided that $S(\zeta)>C_2$, where $C_2$ is another constant).
Lemma \ref{L2} is proved.

\bigskip 

 Lemma~\ref{L2} means that conditions of  Theorem~2.2.7 in~\cite{FMM} are satisfied with the linear function 
$f(\zeta)=S(\zeta)$ and set $A=\{\zeta\in \Z_{+}^{\L}:  S(\zeta)\geq C_1\}$ and  the embedded  Markov chain $\zeta(t)$ is transient in the case $\alpha+\beta\nu> 0$.
\begin{lemma}
\label{L4}
Let $\L$ be a finite connected graph with the constant vertex degree  $\nu(x)\equiv \nu$.
If $\alpha+\beta\nu=0$, then  there exist  $\eps>0$ and $C>0$  such that  the following bound holds
\begin{equation}
\label{2step}
\E(S(\zeta(t+k(\zeta(t)))-S(\zeta(t))\|\zeta(t)=\zeta)\geq \eps,
\end{equation}
provided that $S(\zeta)\geq C$ and 
where 
$$ 
k(\zeta)=
\begin{cases}
1, &  U(x, \zeta)\neq 0\,\,\mbox{for some}\,\, x\in \L, \\
2, & U(x, \zeta)=0\,\,\mbox{for all}\,\, x\in \L.
\end{cases}
$$
\end{lemma}
{\it Proof of Lemma  \ref{L4}.}
As we already noted in the proof of Lemma~\ref{L2}  inequality~(\ref{2step}) is equivalent to the following one 
\begin{equation}
\label{tr2}
{\rm J}(\zeta, \eps)=\delta(\eps)\sum\limits_{x\in \L}e^{U(x, \zeta)} -\sum\limits_{x\in \L}1_{\{\zeta_x>0\}}\geq 0,
\end{equation}
where  $\delta(\eps)$ is defined by~(\ref{delta}) and (\ref{tr2})  would be implied by 
$$\delta(\eps)\sum\limits_{x\in \L}e^{U(x, \zeta)} -|\L|\geq 0.$$
Notice  that by inequality between geometric and arithmetic means we  have  that for all $\zeta$
\begin{equation}
\label{geq}
\sum\limits_{x\in \L}e^{U(x, \zeta)}- |\L|\geq 0,
\end{equation}
since by equation~(\ref{Usum})
\begin{equation}
\label{Usum1}
\sum\limits_{x\in \L} U(x, \zeta)=(\alpha+\beta\nu)S(\zeta)=0.
\end{equation}
It is well known that given numbers $a_1, \ldots, a_m$  geometric and arithmetic means of  these numbers
are equal to each other  if and only if $a_1=\ldots=a_m$.
Therefore, equation~(\ref{Usum1}) also implies that identity 
$ \sum_{x\in \L}e^{U(x, \zeta)}- |\L|= 0$
holds if and only if $U(x, \zeta)=0,$ for all $x\in \L$ otherwise we have got a strict inequality in~(\ref{geq}).
Thus, if there are exactly $0<m\leq |\L|$  vertices with non zero potentials then 
$$\sum\limits_{x\in \L} e^{U(x, \zeta)}- |\L| \geq\sum\limits_{x\in \L:U(x, \zeta) \neq 0}e^{U(x, \zeta)} - m>0.$$
It is easy to see that since the inequality in the preceding display is strict  there exists $\delta_m\in (0, 1)$ such that 
$$\delta_m\sum\limits_{x\in \L} e^{U(x, \zeta)}- |\L| \geq\delta_m\sum\limits_{x\in \L:U(x, \zeta)\neq \neq 0}e^{U(x, \zeta)} - m>0,$$
because  values of potentials $U$ belong to a discrete set $\{ \alpha(k- j/\nu),\, k,j\in \Z_{+}\}$ (where we used that $\beta=-\alpha/\nu$)
which is bounded away from zero.
Thus, given $0<m\leq |\L|$ we claim existence of $\delta_m$ and,   hence, existence  of the corresponding  $\eps=\eps(\delta_m)$
 (using equation~(\ref{delta})).
The required in Lemma~\ref{L4} $\eps$ is obtained  as $\eps=\min_m\eps_m$.

It is easy to see that all potentials cannot stay zero for two steps in a row, hence 
$$
\E\left(S(\zeta(t+2))-S(\zeta(t))\|\zeta(t)=\zeta\right)=\E(S(\zeta(t+2))-S(\zeta(t+1))|\zeta(t)=\zeta) \geq \eps.
$$
Thus inequality~(\ref{2step}) is proven, and by Theorem~2.2.7 in~\cite{FMM}  the embedded Markov chain is transient.  This completes the proof of Lemma~\ref{L4}.

\medskip 

We are ready now to finish the proof of Theorem~\ref{T3}.

\medskip 

{\it Proof of   Part~2) of Theorem~\ref{T3}.} 
 If $\alpha+\beta\nu=0$ then transience of DTMC $\zeta(t)$  implies at least transience of  CTMC $\xi(t)$.
By Theorem~\ref{L1}  CTMC $\xi(t)$ does not explode if  $\alpha<0,\, \alpha+\beta\nu=0$.
Hence,  CTMC $\xi(t)$ is transient if  $\alpha<0,\, \alpha+\beta\nu=0$.
\begin{rem}
{\rm  Let us notice how the sign of parameter $\alpha$ influences  the process dynamics in the case $\alpha+\beta\nu=0$.
 If $\alpha>0,\, \alpha+\beta\nu=0$, then Theorem~\ref{T2} applies  (since $\beta<0$) and, eventually, a single component 
of the Markov chain explodes. A set of configurations  $\{\xi: \xi_x=\xi_y,\, x,y\in\L\}$ is  ''unstable" in the sense that the process tends to leave 
it and to never  return.  In contrast, if $\alpha<0$, then 
the process tends to stay in a neighbourhood  of the same set of configurations (with equal components)
while escaping  to infinity. It is easy to see that vertex potentials are bounded around this set of configurations and this is why 
no explosion happens in this case.
}
\end{rem}
{\it Proof of Part 3) of Theorem~\ref{T3}.}  If  $\alpha+\beta\nu>0$, then explosiveness of CTMC $\xi(t)$ is implied  (regardless of the sign of $\alpha$) 
by Lemma \ref{L3}.  Indeed, given a configuration $\xi$ bound~(\ref{MaxU}) implies the following lower bound for the total transition rate 
$$\sum\limits_{x\in \L}\left(e^{U(x, \xi)}+1_{\{\xi_x>0\}}\right)\geq \max_{x\in \L} e^{U(x,\zeta)}\geq e^{C_1 S(\zeta)},$$
where, as before,  $C_1=(\alpha+\beta\nu)/|\L|$. Besides,  none of the components 
decrease after  time $\tau$  defined in Lemma~\ref{L3}.  Therefore 
the only changes in the systems are jumps up  and these jumps happen with exponentially  increasing rates
whose  inverses   are  summable. 
This yields explosion.

\medskip 

{\it Proof of Part 4) of Theorem~\ref{T3}.}
If both $\alpha>0$ and $\beta>0$ then  transience of DTMC $\zeta(t)$ and explosiveness of CTMC $\xi(t)$ are obvious.
On the other hand, if  $\alpha>\max\{0,\beta\}$ then  Theorem~\ref{T2} applies; if  $0<\alpha<\beta$ and the graph $\L$ is without triangles then Theorem~\ref{T.no.triangle} applies. The proof of Theorem~\ref{T3} is finished.

\subsection{Proof of Theorem~\ref{T31}}

Let $\zeta(t)=(\zeta_1(t), \ldots, \zeta_{n}(t))$ be DTMC corresponding to a complete graph with $n$ vertices. It is easy to see that the potential of a component at vertex $i$ at time $t$ is equal to  
$$
U(i, \zeta(t))=\alpha \zeta_i(t)+\beta(S(\zeta(t))-\zeta_i(t))=(\alpha-\beta)\zeta_i(t)+\beta S(\zeta(t)).
$$
First, we present an intuitive argument justifying the theorem, which is made rigorous later. 
In both cases described in  the theorem, $\alpha+\beta\nu>0$ hence by Lemma~\ref{L3} there exists a moment of time~$\tau$ after which none of the components decrease.  
Also, it is easy to see that in both cases of the theorem $\beta$ must be positive. So, for $t>\tau$ the probability that it is the $i-$th component that increases is equal to
\begin{align*}
\frac{e^{(\alpha-\beta)\zeta_i(t)+\beta S(\zeta(t))}}{\sum_{k=1}^{n} e^{(\alpha-\beta)\zeta_i(t)+\beta S(\zeta(t))}}
&= \frac{e^{(\alpha-\beta)\zeta_i(t)}}{\sum_{k=1}^{n} e^{(\alpha-\beta)\zeta_k(t)}}
\end{align*}
Therefore, in the long run DTMC evolves as a generalized P\'olya urn model with weight function $g(x)=e^{(\alpha-\beta)x}$. Now the well-known results for a generalized P\'olya urn scheme and Theorem~1 in~\cite{Volk1} implies Parts~1) and~3) of Theorem~\ref{T31}. Finally, the explosiveness of the process $\xi(t)$ follows from Parts~3) and~4) of Theorem~\ref{T3}. (One can compare this and the following calculations with the argument presented in the proof of Part~3) of Theorem~\ref{T4}.)

The problem with the above argument is that, strictly speaking, the events $\zeta_{i+1}(t)=\zeta_i(t)+1$, $i=1,2,\dots,n$, and $\tau<t$ are {\em not} independent, as the behaviour of the P\'olya urn {\em may} affect the probability of decreasing of a component. Thus, to make the argument rigorous, we construct the following coupling. 

Let $Y_t$,  $t=1,2,\dots$, be a sequence of i.i.d.\ uniform $[0,1]$ random variables. At time $t$ split the interval $[0,1]$ into $2n$ intervals with lengths proportional to 
$$
[e^{U(1,\zeta(t))},e^{U(2,\zeta(t))},\dots, e^{U(n,\zeta(t))},1,1,\dots,1]
$$
where $U$ is defined by~\eqref{U}. If  $Y_t$ falls into the $i-$th subinterval with $1\le i\le n$ then we set $\zeta_i(t+1)=\zeta_i(t)+1$; if  $n+1\le i\le 2n$ then we set $\zeta_i(t+1)=\max\{0,\zeta_i(t)-1\}$. In both cases we leave the remaining components unchanged. It is easy to see that the process $\zeta(t)$, $t\ge 1$, has exactly the same distribution as the DTMC defined above. At the same time for a fixed $N\in\Z^+$ define the process $\zeta^{(N)}(t)$, $t=N,N+1,\dots$, such that $\zeta^{(N)}(N):=\zeta(N)$ and the transition rules of $\zeta^{(N)}(t)$ are exactly the same as that of $\zeta(t)$ with the only exception that when $Y_t$ falls in the interval with index $\ge n+1$ the process $\zeta^{(N)}(t)$ remains unchanged (i.e., ``no deaths'').
Let $B_N$ be the event ``none of $Y_t$ falls in the intervals indexed $n+1,n+2,\dots,2n$ for all $t\ge N$'', then on $B_N$ we have $\zeta^{(N)}(t)\equiv \zeta(t)$, $t\ge N$, consequently $\zeta(t)$ has the behaviour of the above  P\'olya urn with weight function $g$. Let $A$ be the event $\{\lim_{t\to\infty} \zeta_k(t)/t= 1/n\}$. Since $\zeta_k^{(N)}(t)/t\to 1/n$ a.s., we have
$$
\P(A)\ge \P(A \| B_N) \P(B_N)=\P(B_N).
$$
On the other hand, Lemma~\ref{L3} implies  that $\P(B_N)\to 1$ as $N\to\infty$, which finishes the proof.

\subsection{Proof of Theorem~\ref{T4}}
{\it  Proof of Part~1) of Theorem~\ref{T4}.}  
Throughout the proof, denote the center of the star graph by $n+1$ and all other vertices $1,2,\dots, n$. We skip  the trivial case, where  $\alpha<0$ and $\beta\leq 0$. 

We will show that if  
\begin{equation}
\label{strict1}
\alpha<0<\beta,\,\, \mbox{and}\,\, \alpha+\beta\sqrt{n}<0,
\end{equation}
then the stationary distribution is well defined. 
Let $A_Q$ be the matrix determined by equation~(\ref{AQ}) in the case of the  star graph with $n+1$ vertices.
Denote by $D_n(\mu)$ be the characteristic polynomial  of matrix $A_Q$. A direct computation gives  the  following recursive equation 
$$
D_n(\mu)=(-\alpha-\mu)D_{n-1}(\mu)-\beta^2(-\alpha-\mu)^{n-1},\, \quad n\geq 1,
$$ 
which yields that 
$$
D_n(\mu)=(-1)^{n+1}(\mu+\alpha)^{n-1}(\mu+\alpha+\beta\sqrt{n})(\mu+\alpha-\beta\sqrt{n}).
$$
Thus,  $-\alpha>0$ is the matrix eigenvalue of multiplicity  $n-1$ and  $-\alpha\pm\beta\sqrt{n}$ are eigenvalues of multiplicity  $1$. The eigenvalue $-\alpha-\beta\sqrt{n}>0$ is the minimal one (since $\beta>0$), hence $A_Q$ is positive definite provided 
conditions~(\ref{strict1}) are satisfied.  Positive definiteness of $A_Q$ implies 
 that  $Z_{\alpha, \beta, \L}<\infty$ (as in the proof of Part 1) of Theorem~\ref{T1}). 
Therefore, the stationary distribution is well defined  and the CTMC $\xi(t)$ is ergodic.  

Let us show that if   $\alpha<0<\beta$ and   $\alpha+\beta\sqrt{n}>0$ then $Z_{\alpha, \beta, \L}=\infty$ and the stationary distribution 
is not defined. Start with noticing that  $\left(1, \ldots, 1, \sqrt{n}\right)\in \Z_{+}^{n+1}$ is  the  eigenvector corresponding to the eigenvalue  $(-\alpha-\beta\sqrt{n})$. Therefore, if $\alpha+\beta\sqrt{n}\geq 0$ then  the function $\exp(-Q(\xi)/2)$ is not summable along the direction of this  eigenvalue and, hence, the CTMC $\xi(t)$ is not ergodic. Indeed, in this case, since $\alpha<0$,
\bn
Z_{\alpha,\beta, \L}&= \sum_{\xi\in \Z^{n+1}_{+}} e^{-Q(\xi)/2-\frac{\alpha}{2}\sum_{i=1}^{n+1}\xi_i}
\ge \sum_{\xi\in\Z^{n+1}_{+}\cap G} e^{-Q(\xi)/2}
\en
where $G=\{\xi:\xi_i=[\beta \xi_{n+1}/|\alpha|],\ i=1,2,\dots,n \}$ and $[x]$ denotes the closest integer to $x\in\R$, so that $|x-[x]|\le 1/2$. Using the expression~(\ref{Q}) for $Q(\xi)$ and the fact that $\beta>\sqrt{n} |\alpha|$ we have 
\bn
Z_{\alpha,\beta, \L}& \ge 
\sum_{\xi\in\Z^{n+1}_{+}\cap G} \exp\left(-\frac{n|\alpha|}{8} +\frac{n\beta^2-\alpha^2}{2|\alpha|}\, \xi_{n+1}^2 
\right)
\\ &
=e^{-\frac{n|\alpha|}{8}}\sum_{k=0}^{\infty} \exp\left(
\frac{n\beta^2-\alpha^2}{2|\alpha|} \,k^2 \right)
 =\infty.
\en

\medskip 

{\it Proof of Part 2) of Theorem~\ref{T4}.}
Observe that
\bn
U(i,\zeta)&=-|\alpha|\zeta_{i} +\beta \zeta_{n+1},\ i=1,2,\dots,n;\\
U(n+1,\zeta)&=-|\alpha|\zeta_{n+1} +\beta\sum_{i=1}^n \zeta_i.
\en
An easy calculation gives the following identity 
\begin{equation}
\label{potentials}
(n\beta+|\alpha|)U(n+1,\zeta)+(\beta+|\alpha|)\sum_{i=1}^n U(i,\zeta) =\left(n\beta^2-\alpha^2\right)S(\zeta)
\end{equation}
where $S$ is defined by~(\ref{S}), valid in the case of any star graph. Thus, if   $\alpha+\sqrt{n}\beta=0$,
\bn
& (n\beta+|\alpha|)U(n+1,\zeta)+(\beta+|\alpha|)\sum_{i=1}^n U(i,\zeta) \\ & =
\beta(1+\sqrt{n})\left(\sqrt{n}U(n+1,\zeta)+\sum_{i=1}^n U(i,\zeta)\right)=0
\en
which is equivalent to
\begin{equation}\label{zero}
 \sqrt{n}U(n+1,\zeta)+\sum_{i=1}^n U(i,\zeta)=0.
\end{equation}
Given $\xi$ denote
$$
m_n=m_n(\xi)=\max\limits_{i=1,...,n}\xi_i,
$$
and 
\begin{align*}
\tau_N=\min\{t: \max(\xi_{n+1}(t), \lfloor \sqrt{n}m_n(\xi(t)) \rfloor =N\}
\end{align*}
where $\lfloor a\rfloor\le a$ denotes the integer part of $a$. It is obvious that the Markov chain is explosive if and only if 
$$
\P\left(\lim\limits_{N\to\infty} \tau_N<\infty\right)>0,
$$
but this  cannot happen.
 Indeed, if 
$\xi_{n+1}\geq \lfloor \sqrt{n}m_n \rfloor $ then 
$$
U_{n+1}=\beta(-\sqrt{n}\xi_{n+1}+)\xi_1+\cdots+\xi_n)\leq \sqrt{n}\beta(-\xi_{n+1}+\sqrt{n}m_n)\leq 0,
$$
and, on the other hand, if  $\xi_{n+1}<\lfloor \sqrt{n}m_n\rfloor$ then 
$$
U_k=\beta(-\sqrt{n}m_n+\xi_{n+1})  =
\beta\left[(-\sqrt{n}m_n+ \lfloor \sqrt{n}m_n \rfloor)
-(\lfloor \sqrt{n}m_n \rfloor -\xi_{n+1})\right]<-\beta 
$$
for all $k$ such that $\xi_k=m_n$. Therefore the waiting time $\tau_{N+1}-\tau_N$ is stochastically larger than a certain exponentially distributed random variable which parameters depend only on $n$ and $\beta$ and hence the limit
$\lim_{N\to\infty} \tau_N$  is infinite with probability $1$. 

Now, let us prove transience of DTMC $\zeta(t)$.  Recall  that  $v=(1, \ldots,1,  \sqrt{n})\in \Z_{+}^{n+1}$
 is  the  eigenvector corresponding to the eigenvalue  $(-\alpha-\beta\sqrt{n})$.  Define a function $f$ as the scalar product (in $\R^{n+1}$)
 of vectors  $\zeta$ and $v$, i.e.
$
f(\zeta)=\zeta_1+\ldots+\zeta_n+\sqrt{n}\zeta_{n+1}.
$
For simplicity, denote $f_t=f(\zeta(t))$. We will show that there exists $\eps>0$ such that for all $\zeta$
\begin{align}\label{eqdeltaf}
\E\left[f_{t+2}-f_t\|\zeta(t)=\zeta\right]\geq \eps.
\end{align}
Since the function $f$ is non-negative and has uniformly bounded jumps (as $|f_{t+1}-f_{t}|\le \sqrt{n}$) transience  of $\zeta(t)$ will follow from Theorem~2.2.7 in~\cite{FMM} with $k(\alpha)\equiv 2$. 

To establish~(\ref{eqdeltaf}), observe that for $\eps\in[0,1)$
\begin{align}\label{eqdeltaf2}
& \E\left[f_{t+1}-f_t\|\zeta(t)=\zeta\right]-\eps
\nonumber
\\
&=\frac{\sum_{i=1}^n e^{U(i, \zeta)}+\sqrt{n}e^{U(n+1, \zeta)}-\sum_{i=1}^{n}1_{\{\zeta_{i}>0\}} -  \sqrt{n} 1_{\{\zeta_{n+1}>0\}}}
{\sum_{i=1}^{n+1} \left[ 
e^{U(i, \zeta)}+1_{\{\zeta_{i}>0\}}\right]}-\epsilon
\nonumber \\ & 
=\frac{H(\zeta,\eps)}
{\sum_{i=1}^{n+1} \left[ 
e^{U(i, \zeta)}+1_{\{\zeta_{i}>0\}}\right]}
\end{align}
where
\begin{align*}
H(\zeta,\eps)&=(1-\eps)\sum_{i=1}^n e^{U(i, \zeta)}+\left(\sqrt{n}-\eps\right)e^{U(n+1, \zeta)}
\\ &
-(1+\eps)\sum_{i=1}^{n}1_{\{\zeta_{i}>0\}} -  \left(\sqrt{n}+\eps\right) 1_{\{\zeta_{n+1}>0\}}.
\end{align*}
From~(\ref{zero}) and the inequality between the arithmetical and geometric means we have 
\bn
\sum_{i=1}^n e^{U(i, \zeta)}
\ge n\left[\prod _{i=1}^n e^{U(i, \zeta)}\right]^{1/n}
=n e^{-\frac{U(n+1, \zeta)}{\sqrt{n}}}
\en
hence
\bn
\frac{H(\zeta,\eps)}{1-\eps}&> \sum_{i=1}^n e^{U(i, \zeta)}+\frac{\left(\sqrt{n}-\sqrt{n}\eps\right)}{1-\eps}e^{U(n+1, \zeta)} -\frac{1+\eps}{1-\eps}\sum_{i=1}^{n}1_{\{\zeta_{i}>0\}}
\\ &
-  \frac{\sqrt{n}+\sqrt{n}\eps}{1+\eps} 1_{\{\zeta_{n+1}>0\}}
\\ &
 =\sum_{i=1}^n e^{U(i, \zeta)}+\sqrt{n} e^{U(n+1, \zeta)}
-\frac{1+\eps}{1-\eps}(n+\sqrt{n})
=:\varphi_{\eps}(u)
\en
where
$$
\varphi_{\eps}(u)=n e^{-u/\sqrt n}+\sqrt{n} e^u-\frac{1+\eps}{1-\eps}(n+\sqrt{n})
$$
and $u=U(n+1,\zeta)\in\R$. 

One can easily check that $\varphi_\eps'(0)=0$ and $\varphi_\eps''(u)=e^{-u/\sqrt{n}}+\sqrt{n}e^u>0$ for all $u$, therefore $\varphi_\eps(\cdot)$  attains its unique minimum at $u=0$. If we set $\eps=0$ we also have $\varphi_0(0)=0$ hence $\varphi_0(u)\ge 0$, $u\in \R$ implying that when $\eps=0$  the LHS of~(\ref{eqdeltaf2}) is always non-negative and $f_t$ is thus a submartingale.

To show that it actually increases on average by at least $\eps>0$ in {\em two} steps, note that 
$|U(n+1,\zeta(t+1))-U(n+1,\zeta(t))|\ge \beta>0$ since $\zeta(t+1)$ differs from $\zeta(t)$ in one of the coordinates, and $|\alpha|>\beta$. Therefore,
$$
\min \left\{\left|U\left(n+1,\zeta(t)\right)\right|,
\left|U\left(n+1,\zeta(t+1)\right)\right|\right\}\ge \frac{\beta}2.
$$
Without loss of generality, assume that it is $u=U(n+1,\zeta(t))$ which has the property $|u|\ge \beta/2$. To guarantee that the LHS~(\ref{eqdeltaf2}) is non-negative for some small $\eps>0$ we will establish that 
\begin{align}\label{eqinf}
\inf_{u:\ |u|\ge \beta/2} \varphi_{\eps}(u) =\min\{\phi_{\eps}(-\beta/2),\phi_{\eps}(\beta/2)\}> 0
\end{align}
where the equality follows from the fact that  $\varphi_{\eps}(u)$ is increasing for $u>0$ and decreasing for $u<0$. However, since $\varphi_{0}(\pm\beta/2)$ is strictly positive, as we established before, and $\varphi_{\eps}(u)$ is continuous in $\eps$, by choosing $\eps>0$ sufficiently small we can ensure~(\ref{eqinf}) and hence~(\ref{eqdeltaf}) and transience.

\medskip

{\it Proof of Part 3) of Theorem~\ref{T4}.}
Recall from~(\ref{potentials}) that 
$$
(n\beta+|\alpha|)U(n+1,\zeta)+(\beta+|\alpha|)\sum_{i=1}^n U(i,\zeta) =\left(n\beta^2-\alpha^2\right)S(\zeta),
$$
where now $n\beta^2-\alpha^2>0$, due to our assumption $\beta>|\alpha|/\sqrt{n}$.  Hence, using  the elementary fact that if $a_1+\dots+a_{n+1}=x$ then $\max_i a_i \ge x/(n+1)$ we get that 
\bn
 \max_{i=1,\dots,n+1} U(i,\zeta(k)) 
 \ge CS(\zeta(k))
\en
and $C>0$ is some constant depending on $n$, $\alpha$ and $\beta$.

At the same time, whenever any of the component of $\zeta$ increases, $S(\zeta(k))$ also increases by $1$. For a positive integer $y$ define $\tau_y=\min\{t:\ S(\zeta(t))\ge y\}$. For each  $y\in\{1,2,\dots\}$ the set of configurations of $\zeta$ where $S(\zeta)<y$ is finite, so with probability one at some point of time $k$ the system will reach the state where $S(\zeta(k))\ge y$, consequently $\tau_y<\infty $ a.s.\ for all $y$. Hence we can define the events $A_y=$``there exists $t\ge \tau_y$ such that some component decreases at time $t$''. Then one can easily obtain the following bound
\bn
\P(A_y)\le 
1-\prod_{k=y}^\infty 
\left(1-\frac{n}{e^{\max_i U(\zeta(k),i)}}\right)
\le
1- \prod_{k=y}^{\infty} 
\left(1-\frac{n}{e^{Ck}}\right)
\sim \frac n{1-e^{-C}}\cdot e^{-Cy}
\en
for large enough $y$. Since $\sum_y e^{-Cy}<\infty$ by Borel-Cantelli lemma there will be a.s.\ a time $y'$ for which no $A_y$ ($y\ge \tau_{y'}$) occurs and thus the only changes in the system are increases of the components; this also implies that for any integer  $k>\tau_{y'}$ we  have 
$\max_i U(i,\zeta(k))\ge C(k-k')$, thus ensuring that the CTMC $\xi(t)$ explodes a.s., since the rates of jumps are bounded below by $e^{C(k-k')}$, the inverses of which are again summable.

Let us now observe the DTMC after time $k'$ thus assuming only increases of the components, i.e.\ $S(\zeta(k+1))-S(\zeta(k))=1$ for all $k\ge k'$. Denote
$$
z(k)=\sum_{i=1}^n \zeta_i(k)=S(\zeta(k))-\zeta_{n+1}(k).
$$
Since the probability that only the component at~$n+1$ increases after time~$k$ equals
$$
\prod_{l=k}^{\infty} \frac{e^{U(n+1,\zeta(k))-|\alpha|(l-k)}}{e^{U(n+1,\zeta(k))-|\alpha|(l-k)}
+\sum_{i=1}^n e^{U(i,\zeta(k))}  }=0
$$
on one hand, and  the probability that the component at $n+1$ never increases after time $k$ is equal to
\bn
&\prod_{l=k}^{\infty} \left(
1- \frac{e^{U(n+1,\zeta(k))}}{e^{U(n+1,\zeta(k))}
+\sum_{i=1}^n e^{U(i,\zeta(l))}  }
\right)
\\ &
\le 
\prod_{l=k}^{\infty} \left(
1- \frac{e^{U(n+1,\zeta(k))}}{e^{U(n+1,\zeta(k))}
+n e^{\max_{i=1,\dots,n} U(i,\zeta(l))}  }
\right)\\
&=\prod_{l=k}^{\infty} \left(
1- \left[
1+n \exp\left\{(|\alpha|+2\beta)\zeta_{n+1}(k)-\beta S(\zeta(l))-|\alpha|\min_{i=1,\dots,n} \zeta_i(l))\right\}  \right]^{-1}
\right)
\\ &
\le \prod_{l=k}^{\infty} C\cdot e^{-\beta l}=0
\en
on the other hand, we conclude that both $\zeta_{n+1}(k)\to\infty$ and $z(k)\to\infty$.

Now consider the process $\zeta(k)$ at those times $k_1<k_2<\dots$ when one of the components in $\{1,2,\dots,n\}$  increases. It is easy to see that $z(k_{i+1})-z(k_i)=1$ for all $i$ and that one can couple the process $(\zeta_1(k_i),\zeta_2(k_i),\dots, \zeta_n(k_i))$, $i=1,2,\dots$, with the generalized P\'olya urn with $n$ types of balls and the weight function $g(x)=e^{\alpha x}$. Since $\alpha<0$, from, for example, a trivial comparison with the Friedman urn, we conclude that all $\zeta_j(k_i)$, $j=1,\dots,n$ grow at the same speed, resulting in $\zeta_j(k)/z(k)\to 1/n$. Therefore, for any $\epsilon>0$ there is a (random) time $k_1\ge k'$ such that 
$$
 \frac{1-\epsilon}{n} \le \min_{j=1,\dots,n}\frac{\zeta_j(k)}{z(k)}\le
\max_{j=1,\dots,n}\frac{\zeta_j(k)}{z(k)}\le \frac{1+\epsilon}{n}\text{ for all }k\ge k_1.
$$
Once this being the case, the odds that at time $k$ the component at $n+1$ grows (as opposed to a component at $i$, $i\in\{1,\dots,n\}$) lies in the interval
\bn
\left[
\frac{e^{-|\alpha|\zeta_{n+1}+\beta z}}
{n e^{-|\alpha|(1-\eps) \frac{z}{n}+\beta \zeta_{n+1}} },
\frac{e^{-|\alpha|\zeta_{n+1}+\beta z}}
{n e^{-|\alpha|(1+\eps) \frac{z}{n}+\beta \zeta_{n+1}} }
\right]
=
\left[
 e^{z R_{-\epsilon}-L\zeta_{n+1}-\log(n)},
 e^{z R_{+\epsilon}-L\zeta_{n+1}-\log(n)}
\right]
\en
where
\bn
R_{\pm\epsilon}=\beta+\frac{|\alpha|(1\pm\eps)}{n}, \ L=|\alpha|+\beta.
\en
Let $X(k)=z(k) R_{-\epsilon} -\zeta_{n+1}(k) L$, $k=k_1,k_1+1,\dots$. Then $X(k)$ can be coupled with random walk $Y(k)$ on $[\log(np/(1-p)),+\infty)$ with the transitional probabilities
\bn
Y(k+1)=\begin{cases}
Y(k)+R_{-\epsilon},& \text{ with probability } 1-p;\\
\max\left\{Y(k)-L,\log\left(\frac{np}{1-p}\right)\right\},& \text{ with probability } p,
\end{cases} 
\en
in such a way that $X(k)\le Y(k)$. By choosing $p\in(0,1)$ such that $\E(Y(k+1)-Y(k))=R_{-\epsilon}(1-p) Lp<0$ 
(provided $Y(k)\ge L+\log\left(np/(1-p)\right)$) we ensure that $\lim_{k\to\infty} Y(k)/k=0$, implying in turn that
\bn
\limsup_{k\to\infty} \frac{X(k)}{k}=\limsup_{k\to\infty} \frac{z(k) R_{-\epsilon} -\zeta_{n+1}(k) L}{k} \le 0.
\en
By the completely symmetric argument we also obtain
\bn
\liminf_{k\to\infty} \frac{z(k) R_{+\epsilon} -\zeta_{n+1}(k) L}{k} \ge 0.
\en
Now, using the fact that $z(k)+\zeta_{n+1}(k)=k+const$ for large $k$,
\bn
\frac{R_{-\epsilon}}{L+R_{-\epsilon}} \le \liminf_{k\to\infty} \frac{\zeta_{n+1}(k)}{k}\le
  \limsup_{k\to\infty} \frac{\zeta_{n+1}(k)}{k}
  \le \frac{R_{+\epsilon}}{L+R_{+\epsilon}} 
\en
Since $\epsilon>0$ is arbitrary and $R_{+\epsilon}-R_{-\epsilon}\to 0$ as $\epsilon\to 0$, we get
\bn
\lim_{k\to\infty} \frac{\zeta_{n+1}(k)}{k}=
\frac{\beta+|\alpha|/n}{\beta+|\alpha|/n+\beta+|\alpha|}=\frac{n\beta+|\alpha|}{2n\beta+(n+1)|\alpha|}
\en
and, as a consequence,
\bn
\lim_{k\to\infty} \frac{\zeta_{i}(k)}{k}
=\frac{\beta+|\alpha|}{2n\beta+(n+1)|\alpha|}
\text{ for }i=1,2,\dots,n.
\en
Finally, we also conclude that all the components of the CTMC $\xi$ actually explode simultaneously.

\medskip

{\it Proof of Part 4) of Theorem~\ref{T4}}.
The case $i)$ of the theorem is covered by Theorem~\ref{T2}, and the case $ii)$ is covered by Theorem~\ref{T.no.triangle}, since a star graph does not have triangles.

\end{document}